\documentclass[11pt,a4paper]{article}
\usepackage{babel}
\usepackage{graphicx}
\usepackage{latexsym}
\usepackage{amsmath}
\usepackage{makeidx}
\usepackage[utf8]{inputenc}
\usepackage[utf8]{inputenc}
\usepackage{amsfonts}
\usepackage{amssymb}
\usepackage{amsmath}
\usepackage{amsthm}
\usepackage{tocbibind}
\usepackage{tcolorbox}
\usepackage[framemethod=TikZ]{mdframed}
\usepackage{lipsum}
\usepackage[mathscr]{eucal}
\usepackage{graphicx}
\usepackage{pdfpages}
\usepackage{color}
\usepackage{array}
\usepackage{hyperref}
\usepackage{geometry}

\title{{\bf Strong and weak sharp bounds for Neural Network Operators in Sobolev-Orlicz spaces and their quantitative extensions to Orlicz spaces}}
         
\author{ {\bf Danilo Costarelli}\,\thanks{Corresponding author} \hskip1cm {\bf Michele Piconi}\hskip1cm \\
Department of Mathematics and Computer Science \\
University of Perugia\\
1, Via Vanvitelli, 06123 Perugia, Italy
\\
{\small {\tt danilo.costarelli@unipg.it}} - {\small {\tt michele.piconi@unipg.it}} 
 }
\date{}

\newcommand{\mau}{\geq}
\newcommand{\miu}{\leq}

\newcommand{\N}{\mathbb{N}}
\newcommand{\R}{\mathbb{R}}
\newcommand{\Z}{\mathbb{Z}}
\newcommand{\disp}{\displaystyle}
\newcommand{\be}{\begin{equation}}
\newcommand{\ee}{\end{equation}}
\newcommand{\phis}{\phi_{\sigma}}

\newtheorem{definition}{Definition}[section]
\newtheorem{remark}[definition]{Remark}
\newtheorem{theorem}[definition]{Theorem}
\newtheorem{lemma}[definition]{Lemma}
\newtheorem{proposition}[definition]{Proposition}
\newtheorem{example}[definition]{Example}
\newtheorem{corollary}[definition]{Corollary}

\begin{document}

\maketitle 

\begin{abstract}
In this paper, we establish sharp bounds for a family of Kantorovich-type neural network operators within the general frameworks of Sobolev-Orlicz and Orlicz spaces. We establish both strong (in terms of the Luxemburg norm) and weak (in terms of the modular functional) estimates, using different approaches. The strong estimates are derived for spaces generated by $\varphi$-functions that are $N$-functions or satisfy the $\Delta^\prime$-condition. Such estimates also lead to convergence results with respect to the Luxemburg norm in several instances of Orlicz spaces, including the exponential case. Meanwhile, the weak estimates are achieved under less restrictive assumptions on the involved $\varphi$-function.
To obtain these results, we introduce some new tools and techniques in Orlicz spaces. Central to our approach is the Orlicz Minkowski inequality, which allows us to obtain unified strong estimates for the operators. We also present a weak (modular) version of this inequality holding under weaker conditions.
Additionally, we introduce a novel notion of discrete absolute $\varphi$-moments of the hybrid type, and we employ the Hardy-Littlewood maximal operator within Orlicz spaces for the asymptotic analysis. Furthermore, we introduce the new space $\mathcal{W}^{1,\varphi}(I)$, which is embedded in the Sobolev-Orlicz space $W^{1,\varphi}(I)$ and modularly dense in $L^\varphi(I)$. This allows to achieve asymptotic estimates for a wider class of $\varphi$-functions, including those that do not meet the $\Delta_2$-condition.
For the extension to the whole Orlicz-setting, we generalize a Sobolev-Orlicz density result given by H. Musielak using Steklov functions, providing a modular counterpart. Moreover, we explore the relationships between weak and strong Orlicz–Lipschitz classes, corresponding to the above moduli of smoothness, providing qualitative results on the rate of convergence of the operators. Finally, a (Luxemburg norm) inverse approximation theorem in Orlicz spaces has been established, from which we deduce a characterization of the corresponding Lipschitz classes in terms of the order of convergence of the operators. The latter result shows that some of the achieved estimates are sharp.
 
\end{abstract}

\section{Introduction} \label{sec1}

In the early 1900s, the concept of artificial neural networks (NNs) began to take shape, driven by a desire to replicate the intricate workings of the human brain. This led to the development of network architectures reminiscent of biological neural structures, marking a significant advancement in computational modeling.
\vskip0.2cm

Since the 1950s, the study of NNs has strongly expanded, influencing diverse fields such as neuroscience, engineering, biology, and computer science. This interdisciplinary impact is evident in seminal works like \cite{CY,CAEU1}, which highlight decades of exploration and innovation across various domains. Among the various approaches being explored, a remarkable area of focus is the study of neural network (NN) operators. These have attracted considerable attention, as shown by the numerous papers cited, such as \cite{CACH1,CACH2,COSP2,costarelli2018approximation,KOKR1,ZH1,GNECCObook,kadak2022neural,chen2022construction,qian2022rates,CO2}, with some of them being recent.
\vskip0.2cm

The term {\em NN operators} arose in the paper of Cardaliaguet and Euvrard \cite{CAEU1}, where the authors considered for the first time a usual {\em feed-forward neural network with one hidden layer} activated by $\psi:\R \to \R$ (see, e.g., \cite{BAR,GAO}), of the form:
\be \label{classical-NNs}
\sum_{i=1}^n a_i\, \psi(w_i \cdot x - b_i), \quad x \in \R^d, \quad w_i \in \R^d,\, \quad a_i,\, b_i \in \R, \quad n \in \N,
\ee
expressed as a positive linear operator. Indeed, the definition proposed in \cite{CAEU1} (for functions $f:\R \to \R$ of one variable) was the following:
\be \label{CE-operators}
(B^\alpha_n f)(x)\ :=\ {n^{-\alpha} \over C_b}\, \sum_{k=-n^2}^{n^2} f(k/n)\, b(n^{1-\alpha}(x-k/n)), \quad x \in \R, \quad 0<\alpha<1, \quad n \in \N,
\ee
where $b:\R \to \R^+_0$ is a centered bell-shaped function and $C_b>0$ is a normalization coefficient given by the integral of $b$ on $\R$. It is evident that the operators in (\ref{CE-operators}) are special cases of the neural networks (\ref{classical-NNs}) for some specific choices of the involved parameters. For the linear operators $B^\alpha_n$, it has been proved (see \cite{AN1997}) that any continuous and bounded $f$ can be uniformly approximated on the compact sets with a rate of convergence which is ${\cal O}(\omega(f,n^{\alpha-1}))$, as $n \to +\infty$, where $\omega(f,n^{\alpha-1})$ denotes the usual modulus of continuity of $f$.
\vskip0.2cm

 Later, the definition of the operators $B^\alpha_n$ has been modified in order to consider in place of functions $b$ suitable sigmoidal activation functions, in the exact spirit of the original theory of NNs, as discussed by Cybenko in \cite{CY}. For the latter reason, \cite{CACH1,CACH2} considered respectively NN operators activated by the logistic and the ramp function; however the order of approximation achieved in these paper was ${\cal O}(\omega(f,n^{-1/2}))$, again, only for the case of approximation of continuous functions. The definitions of the NN operators mentioned above have been definitively overcome in \cite{COSP1}, where by the new definition of NN operators an order of approximation of ${\cal O}(\omega(f,n^{-1}))$ was established for $f \in C([a,b])$.

  One of the main peculiarity of {\em all} the above quoted operators was that they were {\em pointwise-type} operators, namely, they depend on the pointwise values assumed by $f$; this makes them not really suitable for the approximation of not necessarily continuous functions, such as $f \in L^p([a,b])$, $1 \miu p<+\infty$. 
\vskip0.2cm

 In the latter context, a possible solution for the problem of the approximation of $f \in L^p([a,b])$ was given by the Kantorovich version of the NN operators introduced and studied in \cite{CS14} (in the case of functions of several variables). 

The key of this approach lies in defining coefficients through local integral averages of an integrable function $f$, in fact removing the pointwise dependence by $f$.
Indeed, the Kantorovich neural network (NN) operators (in the univariate case) take the following form
\begin{equation}\label{NNintro}
(K_n f)(x) = \frac{\displaystyle\sum_{k=\lceil na \rceil}^{\lfloor nb \rfloor-1} \left[ n \int_{k/n}^{(k+1)/n} f(u) du \right] \phis( n x - k )}{\displaystyle\sum_{k=\lceil na \rceil}^{\lfloor nb \rfloor-1} \phis(nx-k)}, 
\end{equation}
where $n \in \mathbb{N}$, $f:[a,b] \to \mathbb{R}$, and $x \in [a,b]$ (with $a, b\in \R$). The function $\phis(\cdot)$ denotes a density function generated by a suitable finite linear combination of sigmoidal activation functions, extensively described in subsequent sections. Additionally, the symbols $\lceil \cdot \rceil$ and $\lfloor \cdot \rfloor$ represent the ceiling and integer part functions, respectively.

  Note that, the Kantorovich NN operators are again special cases of the neural networks in (\ref{classical-NNs}), with the only difference given by the presence of the denominator term. However, the denominator in the definition of (\ref{NNintro}) represents a normalization term, which is useful to get the uniform convergence of $K_nf$ to $f$ on the whole compact $[a,b]$, as discussed in \cite{CACO2025}. Indeed, if the denominator is removed from the definition of $K_nf$, the uniform convergence on the whole $[a,b]$ fails,  but it can be established only on compact sets $I_\delta:=[a+\delta,b-\delta] \subset [a,b]$, with $\delta>0$ sufficiently small. The above mentioned uniform convergence result has been exploited in \cite{CS14} also to prove the $L^p$-convergence thought a density approach.
\vskip0.2cm

 From the mathematical point of view, it is well-known that in the literature non-trivial extensions and generalizations of $L^p$-spaces have been considered; among these we can certainly mention the widely studied Orlicz spaces.
\vskip0.2cm

   The primary aim of this paper is to establish sharp bounds both strong, in terms of the Luxemburg norm, and weak, based on the modular functional, in Orlicz spaces for NN Kantorovich operators outlined in (\ref{NNintro}). This will be done in a quantitative form thanks to the use of strong and weak moduli of smoothness in Orlicz spaces. This error analysis provides a full understanding of the rate of convergence of the above operators, with respect to the Luxemburg norm and the modular functional.  To reach the above purpose, quantitative estimates for the order of approximation in the context of Sobolev-Orlicz spaces must be preliminary proved, in order to extend the result to the whole Orlicz spaces.
To this aim, we need to introduce new tools and techniques of proof that will be described in detail below.
\vskip0.2cm 

  The main framework of this study is represented by the celebrated Orlicz and Sobolev-Orlicz spaces, that have been introduced by J. Musielak and W. Orlicz, in \cite{mo59}. Over the subsequent decades, these spaces garnered considerable attention and were extensively investigated by various mathematical communities worldwide throughout the 1970s-1990s.

  In the late 1970s, Sobolev-Orlicz spaces emerged in literature thanks to T.K. Donaldson in 1971 for problems in functional and differential analysis \cite{D71, DTOS}. Another pioneering research on them between 1976 and 1979 was carried out by H. Hudzik (see, e.g., \cite{h1,h2,h4}). In 1987, V.V. Zhikov highlighted the relevance of Sobolev-Orlicz spaces in elasticity theory, emphasizing their significance in solving Dirichlet problems \cite{Zhikov}, while M. Ruzicka in 2000 (\cite{R2000}) considered applications in modeling electrorheological fluids. These developments attracted a large interest in Sobolev-Orlicz spaces, prompting further research into their applications in partial differential equations (see, e.g., \cite{GMS99, MR08, DF16, CGW21}).
\vskip0.2cm

  In recent years, attention has turned towards extending the study of NN operators to Orlicz spaces, in order to develop a theory covering a wide range of functional spaces. This allows for a unified approach to deriving approximation results across several settings, including $L^p$-spaces, Zygmund (or interpolation) spaces \cite{stein69, BR80}, exponential spaces \cite{H03}, and so on.

  The main purpose of the present paper is to carry on the study of NN operators in the context of Orlicz spaces, started in \cite{CVMathNachr}. Herein, we want to study the asymptotic behaviour of such operators, providing both quantitative and qualitative sharp bounds to fully understand their rate of convergence. To this end, we need to work in the aforementioned framework of Sobolev-Orlicz spaces, where we establish asymptotic results.
\vskip0.2cm

  Specifically, in Section \ref{sec2} we recall the definition of NN Kantorovich operators together with some basic properties. In Section \ref{sec3}, we outline the notations and some preliminary notions on Orlicz and Sobolev-Orlicz setting, that will be used along the paper. 
\vskip0.2cm

  Section \ref{sec4} is devoted to establishing asymptotic estimates within the framework of Sobolev-Orlicz spaces, by using different approaches. We begin by establishing a Luxemburg norm-based inequality, assuming that the density function $\phis$ generating the NN Kantorovich operators is compactly supported and using a Orlicz version of the celebrated Minkowski inequality. To relax the above assumption on $\phis$, we provide a modular asymptotic estimate using a classical moment-type approach and the boundedness of the well-known Hardy-Littlewood maximal function in the Orlicz setting. Here, we need to introduce a new notion of absolute moments of hybrid type, which involve both the $\varphi$-function as the density function. 
The new introduced hybrid-type moment generalizes the classical one (see, e.g., \cite{EstimatesDurr}) connecting the considered density function $\phis$ with the $\varphi$-function $\varphi$ generating the corresponding Sobolev-Orlicz spaces.
Also a stronger version of the well-known $\Delta_2$-condition, i.e., the $\Delta^\prime$-condition, is employed. Moreover, to further weaken the assumptions on the involved $\varphi$-function, we provide another modular estimate which holds in a considerably broader functional setting, including also exponential-type spaces.
Here the price-to-pay for achieving such bounds is that we can consider functions only in a proper subspace of the Sobolev-Orlicz space $W^{1,\varphi}(I)$, that in the paper we denote by $\mathcal{W}^{1,\varphi}(I)$.
\vskip0.2cm

  Section \ref{sectionquantitative} is central, since it focuses on quantitative estimates in terms of the $\varphi$-modulus of smoothness in the whole Orlicz setting. To achieve this, we use the previously established asymptotic estimates in Sobolev-Orlicz spaces as well as density results studied by H. Musielak \cite{Musielak87} in such setting. 
More precisely, we refer to a quantitative density theorem showing that a function in $L^\varphi(I)$ can be approximated by suitable Steklov functions (see \cite{costarelli2024convergence}) with a precise rate given in term of moduli of smoothness. Moreover, we also provided a new modular version of such general result employing a weak version of the modulus of smoothness of $L^\varphi(I)$, defined by means of the modular functional in place of the Luxemburg norm. One of the consequence of such result is that the new space $\mathcal{W}^{1,\varphi}(I)$ is modularly dense in $L^\varphi(I)$.
As a consequence of the above density theorems, we actually obtain both strong, as well as, modular convergence theorems on the whole Orlicz setting. In particular, we are able to obtain a Luxemburg norm convergence result even in the delicate case of exponential-type spaces, whose generating $\varphi$-function does not satisfy the $\Delta_2$-condition.
\vskip0.2cm

 Section \ref{s6} concludes the study with a qualitative analysis within the framework of strong and weak Orlicz Lipschitz classes, defined by using the strong and the weak version of the $\varphi$-modulus of smoothness, respectively. Here, we also explore the inclusion properties of these classes, highlighting the relationships through suitable examples.
\vskip0.2cm

  In Section \ref{sec7}, we provide an inverse approximation result proving a preliminary Bernstein-type inequality in Orlicz spaces with respect to the Luxemburg norm, using the Berens–Lorentz lemma and the equivalence between the notion of the $K$-functional in Orlicz spaces and the strong $\varphi$-modulus of smoothness. 
By combining this with the direct qualitative results, we can also obtain a complete characterization of the strong Orlicz–Lipschitz classes in terms of the rate of convergence of the considered operators. The latter result shows that some of the achieved estimates are sharp.
\vskip0.2cm

  Moreover, along the paper several instances of both sigmoidal functions and Orlicz spaces generated by different $\varphi$-functions are discussed and applied to the main results provided here.
By the way, we also highlight that the well-known Rectified Linear Unit functions (ReLU) can be included in the present treatment, since this can be used to generate certain density functions arising from the so-called ramp functions (see \cite{cpNNlp}).
\vskip0.2cm

 Some final remarks, together with open problems have been discussed in Section \ref{sec8}.  Moreover, we point out that, in general, without requiring any growth condition on the $\varphi$-function, Orlicz and Sobolev-Orlicz spaces are more general than the classical $L^p$ and Sobolev spaces, also providing a counterexample. We also point out the practical relevance of approximation results in Orlicz spaces, which offer greater flexibility than classical $L^p$-approximation in the presence of functions/signals, for instance, with localized peaks or with super-polynomial growth. Finally, we also note that some of the introduced results can be useful also as general results in the context of Sobolev-Orlicz or Orlicz spaces, independently to the applications purposes pursuit in this paper.

\section{The Kantorovich NN operators} \label{sec2} 

A measurable function $\sigma: \R \to \R$ is called a \textit{sigmoidal function} if \[\lim_{x \to -\infty}\sigma(x)=0\quad\text{ and }\quad\lim_{x \to +\infty}\sigma(x)=1.\]
From this point onward, we always consider non-decreasing sigmoidal functions $\sigma$, with $\sigma(1)<1$, satisfying the following conditions:
\begin{itemize}
\item[$(S 1)$] $\sigma(x)-1/2$ is an odd function;
\item[$(S 2)$] $\sigma \in C^2(\R)$ is concave for $x \mau 0$;
\item[$(S 3)$] $\sigma(x)={\cal O}(|x|^{-\alpha-1})$ as $x \to -\infty$, for some $\alpha>0$,
\end{itemize}
as outlined in the general theory provided in \cite{COSP1}.
The density function $\phis$ generated by $\sigma$ is defined as follows
\[
\phi_{\sigma}(x)\, :=\, \frac{1}{2}[\sigma(x+1)-\sigma(x-1)], \hskip1cm x \in \R.
\]
Let now $I:=[a,b]\subset\mathbb{R}$, with $a<b$.
\begin{definition}
Given a sigmoidal function $\sigma$ satisfying the above assumptions and $f:I \to \R$ a locally integrable function, we define the Kantorovich neural network (NN) operators as follows
\begin{equation*} \label{NNopS}
(K_n f)(x)\ =\ {\disp \sum_{k=\lceil na \rceil}^{\lfloor nb \rfloor-1} \left[ n \int_{k/n}^{(k+1)/n} f\left(u  \right)\, du \right] \phis\left( n x - k \right) \over \disp \sum_{k=\lceil na \rceil}^{\lfloor nb \rfloor-1} \phis(nx-k)}, \hskip1cm x \in I,
\end{equation*} 
with $n \in \N$ such that $\lceil na \rceil\le \lfloor nb \rfloor-1$, where $\lceil \cdot \rceil$ and $ \lfloor \cdot \rfloor$ denote the ceiling and the integer part of a given number. 
\end{definition}

The definition of $\phis$ was initially introduced in \cite{CACH1,CACH2} for the logistic and hyperbolic tangent functions, and later it has been extended to all sigmoidal functions $\sigma$ satisfying the previous conditions $(S 1)$, $(S 2)$ and $(S 3)$.

In the following lemma, we outline some auxiliary properties satisfied by $\phis$ which will be used to obtain the results provided in the present paper.

\begin{lemma} \label{lemma1}
$($i$)$ $\phi_{\sigma}(x) \mau 0$ for every $x \in \R$, with $\phi_{\sigma}(2) 
> 0$, and moreover $\disp \lim_{x \to \pm \infty} \phi_{\sigma}(x) = 0$;
\vskip0.2cm
\noindent $($ii$)$ The function $\phi_{\sigma}(x)$ is even;
\vskip0.2cm
\noindent $($iii$)$ The function $\phi_{\sigma}(x)$: 
\begin{equation*} \label{a1}
\mbox{is non-decreasing for}\ x < 0\ \mbox{and 
non-increasing for}\ x \mau 0;
\end{equation*}
\vskip0.2cm
\noindent $($iv$)$ Let $\alpha$ be the positive constant of condition $(S 3)$. Then:
\be \label{a2}
\phi_{\sigma}(x) = \mathcal{O}(|x|^{-\alpha-1}), \quad as \quad x \to 
\pm \infty.
\ee
Hence, it turns out that $\phis \in L^1(\R)$;
\vskip0.2cm
\noindent $($v$)$ For every $x \in \R$,
\be \label{ipsing}
\sum_{k \in \Z} \phis(x-k) = 1,
\ee
and
\begin{equation*}
\|\phis\|_1\ =\ \int_\R \phis(x)\, dx\ =\ 1;
\end{equation*}
\vskip0.2cm
\noindent $($vi$)$ Let $x \in I$ and $n \in \N^+$. Then:
\be \label{a3}
\sum_{k=\lceil na \rceil}^{\lfloor nb \rfloor-1} \phi_{\sigma}(nx - k)\ \mau\ \phi_{\sigma}(2)\ >\ 0.
\ee
\end{lemma}
\begin{remark}\rm
Note that the assumption $\sigma(1)<1$ is needed to prove that $\phis(2)>0$ (see, e.g., \cite{cpNNlp}).
Clearly, if $\sigma$ is strictly increasing, such assumption is no longer required.
\end{remark}

Notice that, $K_n$ are well-defined since its denominator is non-zero by (\ref{a3}) of Lemma \ref{lemma1} and \begin{equation*} \label{well-def}
|(K_n f)(x)|\ \miu\ \|f\|_\infty \ <\ +\infty, \quad x \in I, \quad n \in \N,
\end{equation*}
for any function $f\in L^\infty(I)$, i.e., the space of essentially bounded functions on $I$, endowed with its norm $\|\cdot\|_\infty$. Moreover, if 
we recall the useful notion of the \textit{discrete absolute moment of order $\nu \mau 0$} of $\phis$  (\cite{COSP1}), i.e.,
\begin{equation*}
M_{\nu}(\phis)\ :=\ \sup_{u \in \R}\, \sum_{k \in \Z}\phis(u-k)\, |u-k|^{\nu},
\end{equation*}
we can deduce by (\ref{ipsing}) that 
\begin{equation}\label{m0}
M_{0}(\phis)=1
\end{equation}
 and that \begin{equation*} \label{finite-moments}
M_\nu(\phis)\ <\ +\infty, \hskip1cm 0 \miu \nu < \alpha,
\end{equation*}
see, e.g., \cite{EstimatesDurr}.

\section{Orlicz and Sobolev-Orlicz spaces} \label{sec3}

Here, we now recall some basic notions concerning Orlicz and Sobolev-Orlicz spaces.
\vskip0.2cm

The function $\varphi: \R^+_0 \to \R^+_0$ is said to be a $\varphi$-\textit{function} if it satisfies the following assumptions:
\begin{itemize}
	\item[$(\Phi 1)$] $\varphi \left(0\right)=0$, $\varphi \left(u\right)>0$ for every $u>0$;
	\item[$(\Phi 2)$] $\varphi$ is continuous and non decreasing on $\R^+_0$;
	\item[$(\Phi 3)$] $\disp \lim_{u\to +\infty}\varphi(u)\ =\ + \infty$.
\end{itemize}
For a fixed $\varphi$-function $\varphi$, one can consider the modular $I^{\varphi} : M(I)\to \text{$\R^+_0$}$, where $ M(I)$ denotes the set of all measurable functions $f:I \to [0,+\infty]$, of the form
\begin{displaymath}
I^{\varphi} \left[f\right] := \int_I \varphi(\left| f(x) \right|)\ dx,\ \hskip0.5cm f \in M(I).
\end{displaymath}
The \textit{Orlicz space} generated by $\varphi$ is now defined by
\begin{equation*}
\begin{split}
L^{\varphi}(I) :&= \left\{f \in M\left(I\right):\ I^{\varphi} [\lambda f]<+\infty,\ \mbox{for\ some}\ \lambda>0\right\}.
\end{split}
\end{equation*}
A useful vector subspace of $L^\varphi(I)$ is given by its \textit{space of finite elements}, which is defined as \begin{equation*}
\begin{split}
E^{\varphi}(I) :&= \left\{f \in M\left(I\right):\ I^{\varphi} [\lambda f]<+\infty,\ \mbox{for\ every}\ \lambda>0\right\}.
\end{split}
\end{equation*}
We can introduce in $L^{\varphi}(I)$ a notion of convergence, called \textit{modular convergence}, which induces a topology (modular topology) on the space (\cite{Musielak83, BMV2003}). Namely, we will say that a net of functions $(f_w)_{w>0} \subset L^{\varphi}(I)$ is modularly convergent to a function $f \in L^{\varphi}(I)$ if
\begin{equation}\label{mconv}
\lim_{w \to +\infty} I^{\varphi}\left[\lambda(f_w-f)\right]\ =\ 0,
\end{equation} 
for some $\lambda>0$. In this context, i.e., if the modular $I^\varphi$ is convex (and this is true when $\varphi$ is convex), a notion of norm can be also introduced, namely the so-called \textit{Luxemburg norm}, defined by
\begin{equation}\label{luxemburgnorm}
\|f\|_\varphi:=\inf\left\{u>0:I^{\varphi}\left[\frac{f}{u}\right]\le1\right\}.
\end{equation}
In Theorem 1.1 (b) of \cite{BMV2003}, the authors prove that $\|\cdot\|_\varphi$ is actually a norm, satisfying all the typical properties, such as the absolute homogeneity and the triangle inequality. In this regard, we say that $(f_w)_{w>0} \subset L^{\varphi}(I)$ is convergent to a function $f \in L^{\varphi}(I)$ with respect to the Luxemburg norm, if $\|f-f_w\|_\varphi\to 0$ as $w\to+\infty$, or, equivalently, if (\ref{mconv}) holds for every $\lambda>0$.
\vskip0.2cm

  We say that a $\varphi$-function $\varphi$ verifies the $\Delta_2$-\textit{condition}, if there exists a constant $C_\varphi>0$ such that
\begin{equation}\label{delta2}
\varphi(2u)\le C_\varphi\varphi(u),
\end{equation}
for all $u\ge 0$.

  Obviously, the $\varphi$-function $\varphi(u)=u^p$, with $1\le p<+\infty$, satisfies the $\Delta_2$-condition and, in this case, $L^{\varphi}(I)=L^p(I)$, i.e., it coincides with the classical \textit{Lebesgue spaces}. Other useful examples of Orlicz spaces can be found, e.g., in \cite{mo59,Musielak83,BMV2003,harjulehto2019orlicz} and will be shown later. Moreover, we remark that if the $\varphi$-function $\varphi$ verifies the $\Delta_2$-condition, then the modular and the Luxemburg norm convergences coincide. Furthermore, we also have that $E^\varphi(I)=L^\varphi(I)$.
\vskip0.2cm

  Let now $AC^{(n)}(I)$, with $n\in\mathbb{N}$, be the space of all absolutely continuous functions $f:I\to\mathbb{R}$, having absolutely continuous derivatives up to the order $n$.
For any $\varphi$-function $\varphi$, the space
\begin{equation}\label{defSO}
W^{n,\varphi}(I):=\left\{f\in AC^{(n-1)}(I) \ : \ f^{(n)}\in L^\varphi(I)\right\},
\end{equation}
is called the \textit{Sobolev-Orlicz space} generated by $\varphi$.
In this context, the Luxemburg norm is given by 
\begin{equation*}
\|f\|_{n,\varphi}:=\sum_{\nu=0}^{n}\|f^{(\nu)}\|_\varphi,
\end{equation*}
by which $W^{n,\varphi}(I)$ is actually a Banach space. In particular, if $\varphi(u)=u^p$, with $1\le p\le+\infty$, we obtain the well-known usual \textit{Sobolev spaces} $W^{n,p}(I)$ (see, e.g., \cite{alberico2024modulus}).
\vskip0.2cm

  In order to get quantitative estimates to study the rate of convergence of NN operators recalled in Section \ref{sec2}, it is needed to introduce the notion of $\varphi$-\textit{modulus of smoothness}, that arises from the finite differences of order $k\in\mathbb{N}$ of the considered function, that are defined as
\[
\Delta^k_h(f,x):=\sum_{j=0}^k\binom{k}{j}(-1)^{k-j}f(x+jh),\qquad x\in I,
\]
where, if necessary, $f$ is assumed to be $(b-a)$-periodic (where $b-a$ is the amplitude of the interval $I$).
This leads to the definition of the $\varphi$-\textit{modulus of smoothness of order $k$} of $f$ in Orlicz spaces, namely
\begin{equation}\label{modulusofsmoothnss}
\omega_k(f,\delta)_\varphi:=\sup_{|h|\le\delta}\|\Delta^k_h(f,\cdot)\|_\varphi,\qquad\delta>0.
\end{equation}
We remark that if $k=1$, it reduces to the \textit{first order} $\varphi$-\textit{modulus of smoothness}, that can be briefly denoted by $\omega(f,\delta)_\varphi$.


\section{Asymptotic analysis for the NN operators in Sobolev-Orlicz spaces} \label{sec4}

In this section, we aim to provide an asymptotic analysis for the Kantorovich NN operators in the general setting of Sobolev-Orlicz spaces. Although the context is very general and therefore delicate, working in this framework allows for the development of unifying results valid for a wide range of functional spaces, as we will show later through specific examples.
\vskip0.2cm

  From now on, we always assume that $\varphi$ is a convex $\varphi$-function.

\subsection{A Luxemburg norm-based inequality}

We start from the following result where the Orlicz version of the celebrated Minkowski inequality (see, e.g, p. 144 of \cite{Musielak87}) plays a central role. 

However, to use the aforementioned useful inequality, we must consider Sobolev-Orlicz spaces generated by the so-called \textit{$N$-functions}, that represent a particular class of convex $\varphi$-functions $\varphi$ satisfying
\[
\frac{\varphi(u)}{u}\rightarrow0, \text{ as } u\to0^+\qquad\text{ and }\qquad\frac{\varphi(u)}{u}\to+\infty, \text{ as } u\to+\infty.
\]
We can recall what follows.
\begin{proposition}[Minkowski inequality]\label{propmink} Let $\varphi$ be a given $N$-function. Let $I$ and $J$ be two compact sets of $\mathbb{R}$, $F:I\times J\to\R$ be a measurable function, such that $F(\cdot,t)\in L^\varphi(I)$ a.e. $t\in J$, and $\|F(\cdot,t)\|_\varphi$ is integrable on $J$. Moreover, let $g:J\to \R^+_0$ be integrable on $J$. Thus, it turns out that
\begin{equation}\label{minkineq}    \left \|\int_{J}g(t)\left|F(\cdot,t)\right|dt\right\|_\varphi\le 2\int_{J}g(t) \|F(\cdot,t)\|_\varphi dt.
\end{equation}
\end{proposition}
Now, we are able to prove the following.

\begin{theorem}\label{thm1orlicz} Let $\varphi$ be a given $N$-function and $\sigma$ be a sigmoidal function such that the resulting $\phi_\sigma$ has compact support contained in $[-\Upsilon, \Upsilon]\subset\mathbb{R}$, with $\Upsilon>0$. Thus, for every function $f\in W^{1,\varphi}(I)$, it turns out that
\begin{equation*}\label{luxemburgnormineq}
\left\| K_nf-f\right\|_\varphi\le 4(1+\Upsilon)\cdot \frac{\|f^\prime\|_\varphi}{n},
\end{equation*}
for $n\in\mathbb{N}$. 
\end{theorem}

\begin{proof} Let $f\in W^{1,\varphi}(I)$ and $x\in I$ be fixed. Without any loss of generality, we can extend $f$ as a $(b-a)$-periodic function on the whole $\R$. Moreover, according to the definition of $W^{1,\varphi}(I)$ given in (\ref{defSO}), $f\in AC(I)$, hence it can be expanded by the first-order Taylor formula with integral remainder\footnote{(\ref{taylor1}) is actually the Fundamental Theorem of Calculus for Lebesgue integral.}(see, e.g., \cite{DELO1}), that is
\begin{equation}\label{taylor1}
f(u)=f(x)+\int_x^uf^\prime(t)dt,\qquad u\in I.
\end{equation} 
Hence, we obtain that
\begin{equation}\label{resto}
\begin{split}
\left(K_nf\right)(x)&=f(x)+\ {\disp \sum_{k=\lceil na \rceil}^{\lfloor nb \rfloor-1} \left\{ n \int_{k/n}^{(k+1)/n} \left[\int_x^uf^\prime(t)dt\right]\, du \right\} \phis\left( n x - k \right) \over \disp \sum_{k=\lceil na \rceil}^{\lfloor nb \rfloor-1} \phis(nx-k)}
\\
&=:f(x)+R_{1,n}.
\end{split}
\end{equation}
Now, being the support of $\phi_\sigma$ included in the compact set $[-\Upsilon, \Upsilon]$, we can reduce to consider in the above sum only the integers $k$ such that $|nx-k|\le\Upsilon$. For such integers, by the change of variable $z=t-x$, we can estimate the following
\begin{equation*}
\begin{split}
& \left|n \int_{k/n}^{(k+1)/n} \left[\int_x^uf^\prime(t)dt\right]\, du\right|
\\
& =\left|n \int_{k/n}^{(k+1)/n} \left[\int_0^{u-x}f^\prime(z+x)dz\right]\, du\right|
\\
& \le n \int_{k/n}^{(k+1)/n} \left[\int_{|z|\le{\frac{1}{n}+\frac{|nx-k|}{n}}}\left|f^\prime(z+x)\right|dz\right]\, du
\\
& \le n \int_{k/n}^{(k+1)/n} \left[\int_{|z|\le{\frac{1+\Upsilon}{n}}}\left|f^\prime(z+x)\right|dz\right]\, du
\\
&=\int_{|z|\le{\frac{1+\Upsilon}{n}}}\left|f^\prime(z+x)\right|dz
\end{split}
\end{equation*}
Therefore, by using properties (v) and (vi) of Lemma \ref{lemma1} and the above inequality, the remainder term $R_{1,n}$ can be estimated as follows
\begin{align}\label{R1n}
\left|R_{1,n}(x)\right|&\le\ {\disp \sum_{k=\lceil na \rceil}^{\lfloor nb \rfloor-1} \left[\int_{|z|\le{\frac{1+\Upsilon}{n}}}\left|f^\prime(z+x)\right|dz\right]\, \phis\left( n x - k \right) \over \disp \sum_{k=\lceil na \rceil}^{\lfloor nb \rfloor-1} \phis(nx-k)}\nonumber
\\
&=  \int_{|z|\le{\frac{1+\Upsilon}{n}}}\left|f^\prime(z+x)\right|dz.
\end{align}
Now, we use (\ref{resto}) and the Minkowski inequality in Orlicz spaces (\ref{minkineq}) in order to get
\begin{equation*}
\begin{split}
\|K_n f-f\|_\varphi=\|R_{1,n}(\cdot)\|_\varphi&\le2 \int_{|z|\le{\frac{1+\Upsilon}{n}}}\|f^\prime(z+\cdot)\|_\varphi dz
\\
&=2 \int_{|z|\le{\frac{1+\Upsilon}{n}}}\|f^\prime(\cdot)\|_\varphi dz
\\
&\le 4(1+\Upsilon)\cdot \frac{\|f^\prime\|_\varphi}{n},
\end{split}
\end{equation*}
where $\|f^\prime(\cdot+z)\|_\varphi=\|f^\prime(\cdot)\|_\varphi<+\infty$, being $f^\prime\in L^{\varphi}(I)$ and a $(b-a)$-periodic function. This completes the proof.
\end{proof}

As it is clear, working with the Luxemburg norm is complex, in view of its definition which is based on a computation of an infimum. Based on this motivation, we aim to derive weaker asymptotic estimates by substituting the Luxemburg norm $\|\cdot\|_\varphi$ in the thesis of Theorem \ref{thm1orlicz} with the modular $I^\varphi$, in both the sides of the inequality.
\vskip0.2cm

  To this end, we need to use a weak version of the previous Minkowski inequality. It can be given as follows.
\begin{proposition}[Weak Minkowski inequality]\label{weakMink} Let $\varphi$
 be a convex $\varphi$-function. Let $I$ and $J$ be two compact sets of $\mathbb{R}$, and $F:I\times J\to\R$ be a measurable function, such that $F(\cdot,t)\in L^\varphi(I)$ uniformly with respect to $t\in J$. Moreover, let $g:J\to \R^+_0$ be integrable on $J$. Thus, there exists $\lambda>0$ such that
\begin{equation}\label{weakminkineq}    I^\varphi\left [\lambda \int_{J}g(t)\left|F(\cdot,t)\right|dt\right]\le \frac{1}{\|g\|_1}\int_{J}g(t) \ I^\varphi\left [\lambda \|g\|_1 F(\cdot,t)\right] dt.
\end{equation}    
\end{proposition}
\begin{proof}
Since $F(\cdot,t)\in L^\varphi(I)$ uniformly with respect to $t\in J$, there exists $\lambda^\star>0$ (independent on $t$) such that $I^\varphi\left[\lambda^\star F(\cdot,t)\right]<+\infty$, for every $t\in J$. Let now consider $\lambda>0$ such that $\lambda\|g\|_1\le \lambda^\star$. Applying the well-known Jensen inequality and Fubini-Tonelli theorem, we get
\begin{equation*}
    \begin{split}
        \int_I\varphi\left[\int_J \lambda g(t) \left|F(x,t)\right|\ dt\right]\ dx&\le \int_I\frac{1}{\|g\|_1}\left(\int_J g(t) \ \varphi\left[\lambda \|g\|_1 \left|F(x,t)\right|\right]dt\right) \ dx
        \\
        &=\frac{1}{\|g\|_1}\int_J g(t)\ \left(\int_I\varphi\left[\lambda\|g\|_1 \left|F(x,t)\right|\right]dx\right) \ dt,
    \end{split}
\end{equation*}
which is the desired inequality.
\end{proof}
Now, we are able to prove the following modular version of Theorem \ref{thm1orlicz}.
\begin{theorem}\label{weakthm1orlicz} Let $\varphi$ be a convex $\varphi$-function and $\sigma$ be a sigmoidal function such that the resulting $\phi_\sigma$ has compact support contained in $[-\Upsilon, \Upsilon]\subset\mathbb{R}$, with $\Upsilon>0$. Thus, for every function $f\in W^{1,\varphi}(I)$, there exists $\lambda>0$ such that
\begin{equation*}
I^\varphi\left[\lambda \left(K_nf-f\right)\right]\le I^\varphi\left[\frac{2\lambda (1+\Upsilon)}{n} f^\prime\right],
\end{equation*}
for $n\in\mathbb{N}$. 
\end{theorem}
\begin{proof}
    Since $f\in W^{1,\varphi}(I)$, there exists $\bar{\lambda}>0$ such that $I^\varphi\left[\bar{\lambda}f^\prime\right]<+\infty$. Let us consider $\lambda>0$ such that $2\lambda(1+\Upsilon)<\phis(2)\bar{\lambda}$. Proceeding as in the first part of the proof of Theorem \ref{thm1orlicz} and applying the weak Minkowski inequality stated in (\ref{weakminkineq}) to (\ref{R1n}), we get
    \begin{equation*}
        \begin{split}
            I^\varphi\left[\lambda \left(K_nf-f\right)\right]=I^\varphi\left[\lambda R_{1,n}\right]&\le I^\varphi\left[\int_{|z|\le \frac{1+\Upsilon}{n}}\lambda |f^\prime(z+\cdot)|\ dz\right]
            \\
            &\le\frac{n}{2(1+\Upsilon)}\int_{|z|\le \frac{1+\Upsilon}{n}}I^\varphi\left[\frac{2\lambda(1+\Upsilon)}{n}f^\prime(z+\cdot)\right]\ dz
            \\
            &\le  I^\varphi\left[\frac{2\lambda(1+\Upsilon)}{n}f^\prime\right],
        \end{split}
    \end{equation*}
    being $I^\varphi\left[\bar{\lambda}f^\prime(z+\cdot)\right]=I^\varphi\left[\bar{\lambda}f^\prime(\cdot)\right]<+\infty$. This completes the proof.
\end{proof}
\begin{remark}\rm
Note that, in Theorem \ref{thm1orlicz}, we obtain a strong asymptotic estimate that leads to a Luxemburg norm convergence theorem with order. However, this requires $\varphi$ to be an $N$-function, which excludes some important cases, such as the case $\varphi(u) = u$, $u \ge 0$, namely, the space $W^{1,1}(I)$. Since in Theorem \ref{weakthm1orlicz} $\varphi$ is not necessarily an $N$-function, we can deduce a modular convergence theorem for the cases that are not included in Theorem \ref{thm1orlicz}.\\ Obviously, if $\varphi$ satisfies the $\Delta_2$-condition, also Theorem \ref{weakthm1orlicz} reduces to a Luxemburg norm convergence result. 
\end{remark}


\subsection{A modular asymptotic estimate under a stronger version of the $\Delta_2$-condition}\label{2ndapproach}

 Here, we seek to relax the requirement for the density function $\phis$ to have compact support. This can be achieved by assuming certain moment-type conditions that simultaneously involve the sigmoidal function $\sigma$ and the considered $\varphi$-function $\varphi$.
\vskip0.2cm

   In particular, we introduce a new notion of the absolute moments, which are of the hybrid type.\\ Therefore, we define the \textit{discrete absolute $\varphi$-moment of order $(\mu,\nu)$} as follows
\be\label{phi-moment}
M_{\nu,\mu}^\varphi(\phis)\ :=\ \sup_{u \in \R}\, \sum_{k \in \Z}\phis(u-k)\, |u-k|^{\nu}\varphi\left(|u-k|^\mu\right),\qquad \nu,\mu\ge0.
\ee 
Moreover, we have already discussed about $\Delta_2$-condition (see (\ref{delta2})). Here, we need to use a stronger version of  such condition, which we call as \textit{strong $\Delta_2$-condition} and denote, as in literature, by $\Delta^\prime$ (see, e.g., Definition 7 on p. 28 of \cite{rao1991theory} and also \cite{krivoshein2022wavelet}).\\
We say that a $\varphi$-function $\varphi$ satisfies the \textit{$\Delta^\prime$-condition}, if there exists $C^\prime_\varphi>0$ such that
\begin{equation}\label{delta2strong}
\varphi(uv)\le \text{$C^\prime_\varphi$}\varphi(u)\varphi(v),
\end{equation}
for every $u\ge 0$ and $v\ge 1$.
\vskip0.2cm

  Obviously, taking $v=2$ one can see that the $\Delta^\prime$ implies the $\Delta_2$-condition, with $C_\varphi=C^\prime_\varphi \cdot \varphi(2)$.
\vskip0.2cm

  In the following, we establish an asymptotic modular estimate in Sobolev-Orlicz spaces, which is preparatory for one of the main results presented in the next section. 
\vskip0.2cm

   In order to establish the desired asymptotic formula, we employ the properties of the celebrated \textit{Hardy–Littlewood maximal function} (HL-maximal function), defined as 
\begin{equation*}
	\mathcal{M}f(x):=\sup_{\stackrel{u\in I}{u \neq x} }\frac{1}{\left|x-u\right|}\int_{x}^u\left|f(t)\right|dt,\qquad x\in I,
\end{equation*}
for a locally integrable function $f:I\to\mathbb{R}$ (see, e.g., \cite{stein1970singular, Grafakos2014, MaximalOp}).
\\ In the case of $L^p$-setting, the celebrated theorem of Hardy, Littlewood and Wiener asserts that $\mathcal{M}f$ is bounded on $L^p(\mathbb{R})$ when $1<p\le+\infty$, for every $f\in L^p(\mathbb{R})$, namely
\begin{equation}\label{HLM}
\|\mathcal{M}f\|_p\le C_p\|f\|_p,
\end{equation}
where the constant $C_p>0$ depends only on $p$ (Theorem I.1 of \cite{stein1970singular}).
Extending the above HL-maximal inequality from the $L^p$-setting to the broader framework of Orlicz spaces is not trivial and requires specific conditions on the involved $\varphi$-function.
\vskip0.2cm

   We premise that two $\varphi$-functions $\varphi$ and $\psi$ are \textit{equivalent} if there exists a constant $L\ge 1$ such that 
\begin{equation}\label{equivalence}
\psi\left(\frac{t}{L}\right)\le \varphi(t)\le \psi(Lt),\quad t\ge0.
\end{equation}
  Equivalent $\varphi$-functions lead to equivalent Orlicz spaces (see Theorem 8.17 (a) of \cite{Musielak83}, p. 54).
 \begin{remark}\label{rmkhl} \rm
 In \cite{Hasto15}, the author proves that for any $\varphi$-function $\varphi$ and $\beta>1$ such that 
 \begin{equation}\label{increasing}
     u\mapsto u^{-\beta}\varphi(u), \qquad u>0,
 \end{equation}
 is increasing, there exists another $\varphi$-function $\psi$ equivalent to $\varphi$ (in the sense of (\ref{equivalence})) such that $\psi^{1/\beta}$ is convex. This represents a crucial step to prove the boundedness of the HL-maximal function in Orlicz spaces, which is a highly delicate topic in Functional Analysis (see, e.g., Corollary 3.3 of \cite{Hasto15}).
\end{remark}
\begin{theorem}[\cite{Hasto15}]\label{Mbounded} Let $\varphi$ be a given convex $\varphi$-function and $\beta>1$ so that $u\mapsto u^{-\beta}\varphi(u)$, $u>0$, is increasing. Hence, the HL-maximal function
\begin{equation*}
\mathcal{M}: L^\varphi(\mathbb{R})\to  L^\varphi(\mathbb{R})
\end{equation*}
is bounded, i.e., there exists $K_\varphi>0$ such that $\|\mathcal{M}f\|_\varphi\le K_\varphi\|f\|_\varphi$, for every $f\in L^\varphi(\R)$.
\end{theorem}

\begin{remark}\rm Notice that the existence of $\beta>1$ such that $u\mapsto u^{-\beta}\varphi(u)$, $u>0$, is increasing, turns out to be easily ensured whenever $\varphi$ is an $N$-function.
\end{remark}
Now, we can prove the aforementioned asymptotic formula.
\begin{theorem}\label{thm2} Let $\sigma$ be a sigmoidal function and $\varphi$ be a convex $\varphi$-function satisfying the $\Delta^\prime$-condition and such that $u\mapsto u^{-\beta}\varphi(u)$, $u>0$, is increasing, with $\beta>1$. Moreover, let $M_{0,1}^\varphi(\phis)<+\infty$. Thus, for every function $f\in W^{1,\varphi}(I)$ and $\lambda>0$, there holds
\begin{equation}\label{modineq3}
I^\varphi\left[\lambda\left( K_nf-f\right)\right]\le \psi^{\sigma,\varphi}\cdot I^\varphi\left[\frac{2L^2\lambda}{n} f^\prime\right],
\end{equation}
for $n\in\mathbb{N}$, where $\psi^{\sigma,\varphi}>0$ is a suitable constant and, as stated in Remark \ref{rmkhl}, $L\ge 1$ arises from (\ref{equivalence}). 
\end{theorem}
\begin{proof}
Let $f\in W^{1,\varphi}(I)$ and $n\in\N$ be fixed. 
Firstly, we expand the function $f$ according to the first-order Taylor formula (\ref{taylor1}) as in Theorem \ref{thm1orlicz}, so that
\[
(K_nf)(x)-f(x)=R_{1,n}(x),
\]
for any $x\in I$, where $R_{1,n}$ is given in (\ref{resto}). Then, by using Jensen inequality in its discrete and continuous form (see \cite{COSP2}), and property (\ref{a3}) of Lemma \ref{lemma1}, we obtain
\begin{equation*}
\begin{split}
I^\varphi\left[\lambda\left( K_nf-f\right)\right]&=I^\varphi\left[\lambda\left( R_{1,n}\right)\right]
\\
&\le \int_a^b\varphi\left(\lambda \ {\disp \sum_{k=\lceil na \rceil}^{\lfloor nb \rfloor-1} \left\{ n \int_{k/n}^{(k+1)/n} \left|\int_x^u\left|f^\prime(t)\right|dt\right|\, du \right\} \phis\left( n x - k \right) \over \disp \sum_{k=\lceil na \rceil}^{\lfloor nb \rfloor-1} \phis(nx-k)}\right)\,dx
\\
&\le \frac{\displaystyle\int_a^b\displaystyle\sum_{k=\lceil na \rceil}^{\lfloor nb \rfloor-1}\varphi\left( \lambda n \displaystyle\int_{k/n}^{(k+1)/n} \left|\displaystyle\int_x^u\left|f^\prime(t)\right|dt\right|\, du\right)\phis(nx-k)\,dx}{\disp \sum_{k=\lceil na \rceil}^{\lfloor nb \rfloor-1} \phis(nx-k)}
\\
&\le \frac{1}{\phis(2)}\int_a^b\sum_{k=\lceil na \rceil}^{\lfloor nb \rfloor-1}\left\{n \displaystyle\int_{k/n}^{(k+1)/n}\varphi\left(\lambda\displaystyle\left|\int_x^u\left| f^\prime(t)\right|dt\right|\right)\, du\right\}\phis(nx-k)\,dx
\\
&= \frac{1}{\phis(2)}\int_a^b\sum_{k=\lceil na \rceil}^{\lfloor nb \rfloor-1}\left\{n \displaystyle\int_{k/n}^{(k+1)/n}\varphi\left(\lambda|u-x|\displaystyle\frac{\displaystyle\left|\int_x^u\left| f^\prime(t)\right|dt\right|}{|u-x|}\right)\, du\right\}\phis(nx-k)\,dx
\\
&\le \frac{1}{\phis(2)}\int_a^b\sum_{k=\lceil na \rceil}^{\lfloor nb \rfloor-1}\left\{n \displaystyle\int_{k/n}^{(k+1)/n}\varphi\left(\lambda|u-x|\displaystyle\left|\mathcal{M}f^\prime(x)\right|\right)\, du\right\}\phis(nx-k)\,dx.
\end{split}
\end{equation*}
We now observe that for any $x\in I$ and $u\in\left[\frac{k}{n},\frac{k+1}{n}\right]$, the following inequality holds
\begin{equation}\label{ineq}
\left|u-x\right|\le\left|u-\frac{k}{n}\right|+\left|\frac{k}{n}-x\right|\le\frac{1}{n}+\left|\frac{k}{n}-x\right|.
\end{equation}
Hence, employing (\ref{ineq}), we get
\begin{equation*}
\begin{split}
&I^\varphi\left[\lambda\left( K_nf-f\right)\right]
\\
&\le \frac{1}{\phis(2)}\int_a^b\sum_{k=\lceil na \rceil}^{\lfloor nb \rfloor-1}\left\{n \displaystyle\int_{k/n}^{(k+1)/n}\varphi\left[\lambda\left(\frac{1}{n}+\left|\frac{k}{n}-x\right|\right)\left|\mathcal{M}f^\prime(x)\right|\right]\, du\right\}\phis(nx-k)\,dx
\\
&= \frac{1}{\phis(2)}\int_a^b\sum_{k=\lceil na \rceil}^{\lfloor nb \rfloor-1}\varphi\left[\frac{\lambda}{n}\left(1+\left|k-nx\right|\right)\left|\mathcal{M}f^\prime(x)\right|\right]\phis(nx-k)\,dx.
\end{split}
\end{equation*}
Applying $\Delta^\prime$-condition (see (\ref{delta2strong})) twice, and using the convexity of $\varphi$, there exist $C^\prime_\varphi>0$ such that
\begin{equation*}
\begin{split}
&\varphi\left(\frac{\lambda}{n}\left(1+\left|k-nx\right|\right)\left|\mathcal{M}f^\prime(x)\right|\right)
\\
&\le C^\prime_\varphi\cdot \varphi\left(1+|k-nx|\right)\cdot\varphi\left(\frac{\lambda}{n}\left|\mathcal{M}f^\prime(x)\right|\right)
\\
&\le C^\prime_\varphi\cdot \left[\varphi(2)+\varphi(2|k-nx|)\right]\cdot\varphi\left(\frac{\lambda}{n}\left|\mathcal{M}f^\prime(x)\right|\right)
\\
&\le \varphi(2)\cdot C^\prime_\varphi\cdot \left[1+C^\prime_\varphi\cdot\varphi(|k-nx|)\right]\cdot\varphi\left(\frac{\lambda}{n}\left|\mathcal{M}f^\prime(x)\right|\right).
\end{split}
\end{equation*}

Therefore, recalling (\ref{m0}), we can write
\begin{equation*}
\begin{split}
I^\varphi\left[\lambda\left( K_nf-f\right)\right]
&\le \frac{\varphi(2)\cdot C^\prime_\varphi}{\phis(2)}\int_a^b\sum_{k=\lceil na \rceil}^{\lfloor nb \rfloor-1}\left[1+C^\prime_\varphi\cdot\varphi(|k-nx|)\right]\varphi\left(\frac{\lambda}{n}\left|\mathcal{M}f^\prime(x)\right|\right)\phis(nx-k)\,dx
\\
&\le \frac{\varphi(2)\cdot C^\prime_\varphi}{\phis(2)}\left[1+C^\prime_\varphi\cdot M_{0,1}^\varphi(\phis)\right]\cdot\int_a^b\varphi\left(\frac{\lambda}{n}\left|\mathcal{M}f^\prime(x)\right|\right)\,dx,
\end{split}
\end{equation*}
where $M_{0,1}^\varphi(\phis)<+\infty$ by the assumption.

Now, in order to achieve the thesis, we need to use the auxiliary $\varphi$-function $\psi$ arising from (\ref{equivalence}), such that $\psi^{\frac{1}{\beta}}(x):=\left(\psi(x)\right)^{\frac{1}{\beta}}$ is convex, with $\beta>1$ (see Remark \ref{rmkhl}). Therefore, we can write what follows
 \begin{equation*}
\begin{split}
I^\varphi\left[\lambda\left( K_nf-f\right)\right]&\le  \frac{\varphi(2)\cdot C^\prime_\varphi}{\phis(2)}\left[1+C^\prime_\varphi\cdot M_{0,1}^\varphi(\phis)\right]\int_a^b\varphi\left[\frac{\lambda}{n}\left|\mathcal{M}f^\prime(x)\right|\right]\,dx
\\
&\le \frac{\varphi(2)\cdot C^\prime_\varphi}{\phis(2)}\left[1+C^\prime_\varphi\cdot M_{0,1}^\varphi(\phis)\right]\int_a^b\psi\left[\frac{L\lambda}{n}\left|\mathcal{M}f^\prime(x)\right|\right]\,dx
\\
&= \frac{\varphi(2)\cdot C^\prime_\varphi}{\phis(2)}\left[1+C^\prime_\varphi\cdot M_{0,1}^\varphi(\phis)\right]\int_a^b\left(\psi^{\frac{1}{\beta}}\left[\frac{L\lambda}{n}\left|\mathcal{M}f^\prime(x)\right|\right]\right)^\beta\,dx.
\end{split}
\end{equation*}
Now, we can apply Jensen inequality in its continuous form to the convex function $\psi^{\frac{1}{\beta}}$ as follows
\begin{equation*}
\begin{split}
\psi^{\frac{1}{\beta}}\left[\frac{L\lambda}{n}\left|\mathcal{M}f^\prime(x)\right|\right]
&=\psi^{\frac{1}{\beta}}\left[\frac{L\lambda}{n}\sup_{\stackrel{u\in I}{u \neq x}}\frac{\left|\displaystyle\int_x^u\left|f^\prime(t)\right| \ dt\right|}{|u-x|}\right]
\\
&\le \sup_{\stackrel{u\in I}{u \neq x}}\frac{\left|\displaystyle\int_x^u\psi^{1/\beta}\left[\frac{2L\lambda}{n}\left|f^\prime(t)\right| \right] \ dt\right|}{|u-x|}
\\
&=\mathcal{M}\left(\psi^{1/\beta}\left[\frac{2L\lambda}{n}\left|f^\prime(x)\right|\right]\right),
\end{split}
\end{equation*}
where in the above inequalities we used the definition of the supremum.
Observing now that $\psi^{1/\beta}\left[\frac{2L\lambda}{n}\left|f^\prime\right|\right]\in L^\beta(I)$, with $\beta>1$, from (\ref{HLM}) there exists $C_\beta>0$ such that
\begin{equation*}
\begin{split}
I^\varphi\left[\lambda\left( K_nf-f\right)\right]
&\le  \frac{\varphi(2)\cdot C^\prime_\varphi}{\phis(2)}\left[1+C^\prime_\varphi\cdot M_{0,1}^\varphi(\phis)\right]\int_a^b\left(\mathcal{M}\left(\psi^{1/\beta}\left[\frac{2L\lambda}{n}\left|f^\prime(x)\right|\right]\right)\right)^\beta\,dx
\\
&\le  C_\beta^\beta\cdot \frac{\varphi(2)\cdot C^\prime_\varphi}{\phis(2)}\left[1+C^\prime_\varphi\cdot M_{0,1}^\varphi(\phis)\right]\int_a^b\psi\left[\frac{2L\lambda}{n}\left|f^\prime(x)\right|\right]\,dx
\end{split}
\end{equation*}
\begin{equation*}
\begin{split}
&\le C_\beta^\beta\cdot \frac{\varphi(2)\cdot C^\prime_\varphi}{\phis(2)}\left[1+C^\prime_\varphi\cdot M_{0,1}^\varphi(\phis)\right]\int_a^b\varphi\left[\frac{2L^2\lambda}{n}\left|f^\prime(x)\right|\right]\,dx
\\
&=:\psi^{\sigma,\varphi}\cdot I^\varphi\left[\frac{2L^2\lambda}{n}f^\prime\right],
 \end{split}
\end{equation*}
in view of (\ref{equivalence}). 

Since the $\Delta^\prime$-condition holds, and consequently the $\Delta_2$-condition is also valid, we can conclude that \(f^\prime \in L^\varphi(I) = E^\varphi(I)\). Therefore, \(\displaystyle I^\varphi\left[\frac{2L^2\lambda}{n}f^\prime\right] < +\infty\) for every \(\lambda > 0\). This completes the proof.
\end{proof}
Note that Theorem \ref{thm2} is actually a Luxemburg norm convergence theorem with order. Indeed, using Theorem \ref{thm2}, we can derive the following.
\begin{theorem}\label{newthm}
    Under the assumptions of Theorem \ref{thm2}, for any $f \in W^{1,\varphi}(I)$ there holds
    \begin{equation*}\label{strongineq2}
\| K_nf-f \|_\varphi\le M\frac{2L^2}{n} \|f^\prime\|_\varphi,
\end{equation*}
where $M=\max\left\{1,\psi^{\sigma,\varphi}\right\}, L\ge 1$ are the constants arising from Theorem \ref{thm2}.
\end{theorem}
\begin{proof}
    Let $f \in W^{1,\varphi}(I)$. By Theorem \ref{thm2}, we know that the modular asymptotic estimate given in (\ref{modineq3}) holds for every $\lambda > 0$ and $n\in\N$. Let now $\lambda^\star:=\frac{1}{\left\|\frac{2L^2}{n}f^\prime\right\|_\varphi}>0$ be such that $I^\varphi\left[\lambda^\star \frac{2L^2}{n} f^\prime\right]\le 1$. Obviously, here we consider the nontrivial case $\|f^\prime\|_\varphi\ne0$. Therefore, from (\ref{modineq3}) we obtain that
\[
I^\varphi\left[\lambda^\star\left(K_nf-f\right)\right]\le \psi^{\sigma,\varphi} \cdot I^\varphi\left[\frac{2L^2\lambda^\star}{n} f^\prime\right]\le \max\left\{1,\psi^{\sigma,\varphi}\right\}=:M,
\]
for $n\in\N$. Now, by the convexity of the modular, we have
\[
I^\varphi\left[\frac{\lambda^\star\left(K_nf-f\right)}{M}\right]\le\frac{1}{M}\cdot I^\varphi\left[\lambda^\star\left(K_nf-f\right)\right]\le 1,
\]
and so by (\ref{luxemburgnorm}), we observe that
\[
\left\|\lambda^\star\left(K_nf-f\right)\right\|_\varphi=\inf\left\{u>0:I^\varphi\left[\frac{\lambda^\star \left(K_nf-f\right)}{ u}\right]\le 1\right\}\le M,
\]
from which we finally obtain
\begin{equation*}
    \| K_nf-f \|_\varphi\le M\frac{2L^2}{n} \|f^\prime\|_\varphi.
\end{equation*}
\end{proof}
\begin{remark}\rm
    We observe that in Theorem \ref{newthm} (respect to Theorem \ref{thm1orlicz}) the requirement for the density function $\phis$ to be compactly supported has been eliminated, assuming that the involved $\varphi$-function satisfies the $\Delta^\prime$-condition and (\ref{increasing}).
\end{remark}


\subsection{A more general modular asymptotic-type theorem}\label{hl}

The aim of this subsection is to establish a modular asymptotic-type theorem in Orlicz spaces generated by $\varphi$-functions that, in general, do not satisfy $\Delta^\prime$, as well as $\Delta_2$-condition.
\vskip0.2cm

   In this context, we need to define a suitable subspace  of $W^{1,\varphi}(I)$, namely 
\begin{equation}\label{subspace}
    \mathcal{W}^{1,\varphi}(I)=\left\{f\in W^{1,\varphi}(I) \ : \ \left(f^\prime\right)^2\in L^\varphi(I)\right\}.
    \end{equation}
In general, it is clear that the following proper inclusion holds
\[
\mathcal{W}^{1,\varphi}(I)\subset W^{1,\varphi}(I).
\]

We can prove the following.

\begin{theorem}\label{thm3} Let $\sigma$ be a sigmoidal function and $\varphi$ be a convex $\varphi$-function such that $u\mapsto u^{-\beta}\varphi(u)$, $u>0$, is increasing, with some $\beta>1$. Moreover, let $M_{0,2}^\varphi(\phis)<+\infty$. Then, for every $f\in \mathcal{W}^{1,\varphi}(I)$, there exists $0<\lambda<1$ such that
\[
I^\varphi\left[\lambda\left(K_nf-f\right)\right]\le \xi_1^{\sigma,\varphi}\cdot\frac{\lambda}{n}+\xi_{2}^{\sigma,\varphi}\cdot I^\varphi\left[\frac{4\lambda L^2}{n} \left(f^\prime\right)^2\right],
\]
for $n\in\N$ sufficiently large, where $\xi_1^{\sigma,\varphi}>0$ and $\xi_{2}^{\sigma,\varphi}>0$ are suitable constants and, as established in Remark \ref{rmkhl}, $L\ge 1$ arises from (\ref{equivalence}).
\end{theorem}
\begin{proof}
Let $f\in \mathcal{W}^{1,\varphi}(I)$ and $n\in\N$ be fixed. By definition, being $\left(f^\prime\right)^2\in L^\varphi(I)$, there exists $\bar{\lambda}>0$ such that
\begin{equation}\label{lambdabar3}
I^\varphi\left[\bar{\lambda} \left(f^\prime\right)^2\right]<+\infty.
\end{equation}
Let us consider $\lambda>0$ such that
\begin{equation}\label{lambda3}
4\lambda L^2<\min\left\{\bar{\lambda},1\right\},
\end{equation}
where $L\ge 1$ arises again from (\ref{equivalence}). 
Proceeding as in the proof of Theorem \ref{thm2}, we arrive to
\begin{equation*}
\begin{split}
I^\varphi\left[\lambda\left( K_nf-f\right)\right]&\le \frac{1}{\phis(2)}\int_a^b\sum_{k=\lceil na \rceil}^{\lfloor nb \rfloor-1}\varphi\left(\frac{\lambda}{n}\left(1+\left|k-nx\right|\right)\left| \mathcal{M}f^\prime(x)\right|\right) \phis(nx-k)\,dx.
\end{split}
\end{equation*}

We now consider the trivial inequality $c d\le\frac{c^2+d^2}{2}$, with $c$, $d>0$. Taking $c:=\sqrt{\frac{\lambda}{n}}\left(1+|k-nx|\right)$ and $d:=\sqrt{\frac{\lambda}{n}}\left|\mathcal{M}f^\prime(x)\right|$ and employing the convexity of $\varphi$, we can write 
\begin{equation*}
\begin{split}
\varphi\left(\sqrt{\frac{\lambda}{n}}\left(1+|k-nx|\right)\cdot \sqrt{\frac{\lambda}{n}}\left|\mathcal{M}f^\prime(x)\right| \right)\le \varphi\left(\frac{\lambda}{n}\left(1+|k-nx|\right)^2\right)+\varphi\left(\frac{\lambda}{n}\left|\mathcal{M}f^\prime(x)\right|^2\right).
\end{split}
\end{equation*}
Hence, exploiting the above inequality, the convexity of $\varphi$ and (\ref{m0}), we obtain
\begin{equation*}
\begin{split}
I^\varphi\left[\lambda\left( K_nf-f\right)\right]&\le \frac{1}{\phis(2)}\int_a^b\sum_{k=\lceil na \rceil}^{\lfloor nb \rfloor-1}\varphi\left(\frac{\lambda}{n}\left(1+|k-nx|\right)^2\right) \phis(nx-k)\,dx
\\
&\qquad+\frac{1}{\phis(2)}\int_a^b\sum_{k=\lceil na \rceil}^{\lfloor nb \rfloor-1}\varphi\left(\frac{\lambda}{n}\left|\mathcal{M}f^\prime(x)\right|^2\right) \phis(nx-k)\,dx
\\
&\le \frac{1}{\phis(2)}\int_a^b\sum_{k=\lceil na \rceil}^{\lfloor nb \rfloor-1}\left[\varphi\left(\frac{2\lambda}{n}\right)+\varphi\left(\frac{2\lambda}{n}\left(k-nx\right)^2\right)\right] \phis(nx-k)\,dx
\\
&\qquad+\frac{1}{\phis(2)}\int_a^b\sum_{k=\lceil na \rceil}^{\lfloor nb \rfloor-1}\varphi\left(\frac{\lambda}{n}\left|\mathcal{M}f^\prime(x)\right|^2\right) \phis(nx-k)\,dx
\\
&\le \frac{1}{\phis(2)}\int_a^b\sum_{k=\lceil na \rceil}^{\lfloor nb \rfloor-1}\frac{2\lambda}{n}\left[\varphi\left(1\right)+\varphi\left(\left(k-nx\right)^2\right)\right] \phis(nx-k)\,dx
\\
&\qquad+\frac{M_0(\phis)}{\phis(2)}\int_a^b\varphi\left(\frac{\lambda}{n}\left|\mathcal{M}f^\prime(x)\right|^2\right) \,dx
\\
&\le \frac{2\lambda}{n}\cdot\frac{(b-a)}{\phis(2)}\left[\varphi\left(1\right)+M_{0,2}^\varphi(\phis)\right]
+\frac{1}{\phis(2)}\int_a^b\varphi\left(\frac{\lambda}{n}\left|\mathcal{M}f^\prime(x)\right|^2\right) \,dx
\\
&=: \frac{2\lambda}{n}\cdot\frac{(b-a)}{\phis(2)}\left[\varphi\left(1\right)+M_{0,2}^\varphi(\phis)\right]+\mathcal{T}
\\
&=:\xi_1^{\varphi,\sigma}\cdot\frac{\lambda}{n}+\mathcal{T},
\end{split}
\end{equation*}
for $n\in\N$. 
\\
Now, we have to focus on the term $\mathcal{T}$. Here, we need to use again the auxiliary $\varphi$-function $\psi$ arising from (\ref{equivalence}) such that $\psi^{1/\beta}$, with $\beta>1$, is convex. Hence, there exists $L\ge 1$ such that
\begin{equation*}
\begin{split}
\mathcal{T}&\le  \frac{1}{\phis(2)}\int_a^b\varphi\left[\frac{\lambda}{n}\left|\mathcal{M}f^\prime(x)\right|^2\right]dx
\\
&\le \frac{1}{\phis(2)}\int_a^b\psi\left[\frac{\lambda L}{n}\left|\mathcal{M}f^\prime(x)\right|^2\right]dx
\\
&= \frac{1}{\phis(2)}\int_a^b\left(\psi^\frac{1}{\beta}\left[\frac{\lambda L}{n}\left|\mathcal{M}f^\prime(x)\right|^2\right]\right)^\beta dx.
\end{split}
\end{equation*}
Now, using Jensen inequality twice in its continuous form to the convex functions $|\cdot|^2$ and $\psi^{\frac{1}{\beta}}$, respectively, we have
\begin{equation*}
\begin{split}
\psi^{\frac{1}{\beta}}\left[\frac{ \lambda L}{n}\left|\mathcal{M}f^\prime(x)\right|^2\right]
&=\psi^{\frac{1}{\beta}}\left[\frac{\lambda L}{n}\left|\sup_{\stackrel{u\in I}{u \neq x}}\frac{\left|\displaystyle\int_x^u\left|f^\prime(t)\right| \ dt\right|}{|u-x|}\right|^2\right]
\\
&\le\psi^{\frac{1}{\beta}}\left[\frac{2\lambda L}{n}\sup_{\stackrel{u\in I}{u \neq x}}\frac{\left|\displaystyle\int_x^u\left|f^\prime(t)\right|^2 \ dt\right|}{|u-x|}\right]
\\
&\le \sup_{\stackrel{u\in I}{u \neq x}}\frac{\left|\displaystyle\int_x^u\psi^{\frac{1}{\beta}}\left[\frac{4\lambda L}{n}\left|f^\prime(t)\right|^2\right] \ dt\right|}{|u-x|}
\\
&=\mathcal{M}\left(\psi^{1/\beta}\left[\frac{4\lambda L}{n}\left|f^\prime(x)\right|^2\right]\right).
\end{split}
\end{equation*}
From (\ref{equivalence}), we know that $\psi^{1/\beta}\left[\frac{4\lambda L}{n}\left|f^\prime\right|^2\right]\in L^\beta(I)$, with $\beta>1$, hence, by (\ref{HLM}) there exists a constant $C_\beta>0$ such that
\begin{equation*}
\begin{split}
\mathcal{T}&\le  \frac{C_\beta^\beta}{\phis(2)}\int_a^b\psi\left[\frac{4\lambda L}{n}\left|f^\prime(x)\right|^2\right] dx
\\
&\le  \frac{C_\beta^\beta}{\phis(2)}\int_a^b\varphi\left[\frac{4\lambda L^2}{n}\left|f^\prime(x)\right|^2\right] dx
\\
&=:\xi_2^{\varphi,\sigma}\cdot I^\varphi\left[\frac{4\lambda L^2}{n}(f^\prime)^2\right],
\end{split}
\end{equation*}
in view of (\ref{equivalence}), (\ref{lambda3}) and (\ref{lambdabar3}). This completes the proof.
\end{proof}

Note that Theorem \ref{thm3} is a modular convergence theorem with order, from which, in general, it is not possible to deduce a corresponding inequality involving the Luxemburg norm.


\section{Quantitative analysis in Orlicz spaces}\label{sectionquantitative}

 Here, to establish estimates based on the $\varphi$-modulus of smoothness defined in (\ref{modulusofsmoothnss}), we employ a result given by H. Musielak in 1987 (see \cite{Musielak87}) by using Steklov functions (see \cite{costarelli2024convergence}). The latter extends a theorem that B. Sendov and V.A. Popov provided in $L^p$-spaces (see \cite{SP88}) to the broader context of Orlicz spaces.

\begin{theorem}[\cite{Musielak87}]\label{SPOrlicz} Let $\varphi$ be an $N$-function and $f$ be a $(b-a)$-periodic function belonging to $L^\varphi(I)$. Let $k>0$ be an integer and $0<h\le (b-a)/k$.
Let $f_{k,h}\in L^\varphi(I)$ be the Steklov function defined by
\begin{equation}\label{steklov}
\begin{split}
f_{k,h}(x)&=(-h)^{-k}\int_0^h\dots\int_0^h \sum_{m=1}^k(-1)^{k-m+1}\binom{k}{m}f\left(x+\frac{m}{k}(t_1+\dots+t_k)\right) dt_1\dots dt_k.
\end{split}
\end{equation}
Thus, it turns out that
\vskip0.2cm
\noindent $($i$)$ $\|f-f_{k,h}\|_\varphi\le 2\omega_k(f,h)_\varphi$;
\vskip0.2cm
\noindent $($ii$)$ $f_{h,k}\in W^{k,\varphi}(I)$ and \[\|f_{k,h}^{(s)}\|_\varphi\le 2(2k)^kh^{-s}\omega_s(f,h)_\varphi,\qquad s=1,2,\dots,k.\]
\end{theorem}

In order to achieve the first of the main results of the present section, we premise the following preliminary theorem, showing that the family of Kantorovich NN operators turns out to be well-defined in Orlicz spaces.

\begin{lemma}[\cite{CVMathNachr}]\label{lemmaorlicz} For any $f\in L^\varphi(I)$ with $\varphi$ a convex $\varphi$-function and $\lambda>0$, there holds
\begin{equation}\label{modularineqlemma}
I^\varphi\left[\lambda K_nf\right]\le\phi(2)^{-1}I^\varphi\left[\lambda f\right].
\end{equation}
\end{lemma}

\begin{remark}\label{rmk1}\rm 
It is important to outline that the modular inequality for Kantorovich NN operators stated in (\ref{modularineqlemma}) leads to the corresponding Luxemburg norm version, given by
\begin{equation}\label{strongineq1}
\| K_nf\|_\varphi\le \phis(2)^{-1}\|f\|_\varphi,
\end{equation}
for any $f\in L^\varphi(I)$. Indeed, let $f\in L^\varphi(I)$ be fixed. Therefore, by the definition of Luxemburg norm given in (\ref{luxemburgnorm}), there exists $\lambda^\star:=\frac{1}{\|f\|_\varphi}>0$ so that $I^\varphi[\lambda^\star f]\leq 1$. Clearly, we here exclude again the trivial case $\|f\|_\varphi=0$.
Thus, by (\ref{modularineqlemma}) we have that
\[
I^\varphi[\lambda^\star K_nf]\le \phi(2)^{-1} I^\varphi[\lambda^\star f]\leq \phi(2)^{-1}.
\]
 Note that, by the properties of $\sigma$, we have that $\phis(2)^{-1}>1$. Hence, by the convexity of the modular $I^\varphi$, we have that
\[
I^\varphi\left[\frac{\lambda^\star K_n f}{\phis(2)^{-1}}\right]\le\frac{I^\varphi\left[\lambda^\star K_n f\right]}{\phis(2)^{-1}}\le 1.
\]
 By using again the definition of Luxemburg norm and the latter inequality, we can state that
\[\|\lambda^\star K_n f\|_\varphi=
\inf\left\{u>0:I^\varphi\left[\frac{\lambda^\star K_n f}{u}\right]\le 1\right\}\le \phis(2)^{-1}.
\]
Now, recalling that $\varphi$ is convex and so $\|\cdot\|_\varphi$ is a norm, we immediately get

\[
\| K_nf\|_\varphi\le\phis(2)^{-1} \|f\|_\varphi.
\]

\end{remark}


Now, we are ready to give the desired quantitative estimate based on the $\varphi$-modulus of smoothness $\omega(f, \cdot)_\varphi$ defined in (\ref{modulusofsmoothnss}).

\begin{theorem}\label{mainthm1orlicz} Let $\varphi$ be an $N$-function and $\sigma$ be a sigmoidal function such that the resulting $\phi_\sigma$ has compact support in $[-\Upsilon, \Upsilon]\subset\mathbb{R}$, with $\Upsilon>0$. Thus, for every function $f\in L^{\varphi}(I)$, there holds
\begin{equation*}
\| K_nf-f\|_\varphi\le\Lambda \cdot \omega\left(f,\frac{1}{n}\right)_\varphi,
\end{equation*}
where $\Lambda>0$ is a suitable constant.
\end{theorem}
\begin{proof}
Let $f\in L^{\varphi}(I)$ be fixed. Without any loss of generality, we can extend $f$ to the whole $\R$ as a $(b-a)$-periodic function. By Theorem \ref{SPOrlicz} (i) and (ii), there exists $f_{1,h}\in W^{1,\varphi}(I)$ such that
\begin{equation}\label{iSPOrlicz}
\|f-f_{1,h}\|_\varphi\le 2 \ \omega(f,h)_\varphi,\qquad 0<h\le(b-a).
\end{equation}
Moreover, Lemma \ref{lemmaorlicz} and (\ref{strongineq1}) lead to
\begin{equation}\label{NNOrlicz}
\| K_n\left(f-f_{1,h}\right)\|_\varphi\le \phis(2)^{-1}\|f-f_{1,h}\|_\varphi,
\end{equation}
 being $f-f_{1,h}\in L^\varphi(I)$.\\
Then, by the linearity of the NN Kantorovich operators, the properties of $\|\cdot\|_\varphi$, and in view of (\ref{iSPOrlicz}) and (\ref{NNOrlicz}), we have
\begin{equation*}
\begin{split}
\| K_nf-f\|_\varphi&\le \| K_nf-K_nf_{1,h}\|_\varphi+\|K_nf_{1,h}-f_{1,h}\|_\varphi+\|f_{1,h}-f\|_\varphi
\\
&\le \left(\phis(2)^{-1}+1\right)\|f_{1,h}-f\|_\varphi+\|K_nf_{1,h}-f_{1,h}\|_\varphi
\\
&\le 2\left(\phis(2)^{-1}+1\right)\omega(f,h)_\varphi +\|K_nf_{1,h}-f_{1,h}\|_\varphi.
\end{split}
\end{equation*}
Now, by using Theorem \ref{thm1orlicz} and, consequently, Theorem \ref{SPOrlicz} (ii), we can write
\begin{equation*}
\begin{split}
\|K_nf_{1,h}-f_{1,h}\|_\varphi&\le 4(1+\Upsilon)\cdot \frac{\|f_{1,h}^\prime\|_\varphi}{n}
\\
 &\le 16(1+\Upsilon)\cdot \frac{h^{-1}}{n}\cdot\omega(f,h)_\varphi.
\end{split}
\end{equation*}
In summary, the following estimate holds
\begin{equation*}
\begin{split}
\| K_nf-f\|_\varphi&\le  2\left(\phis(2)^{-1}+1\right)\omega(f,h)_\varphi+\frac{16(1+\Upsilon) h^{-1}}{n}\cdot \omega(f,h)_\varphi.
\end{split}
\end{equation*}
Now, considering $h=\frac{1}{n}\le (b-a)$ with $n\in\mathbb{N}$, we finally get 
\begin{equation*}
\begin{split}
\| K_nf-f\|_\varphi&\le \left\{ 2\left(\phis(2)^{-1}+1\right) +16(1+\Upsilon) \right\} \cdot \omega\left(f,\frac{1}{n}\right)_\varphi=:\Lambda \cdot \omega\left(f,\frac{1}{n}\right)_\varphi,
\end{split}
\end{equation*}
for every sufficiently large $n\in\mathbb{N}$. Hence, the proof is now completed.
\end{proof}

Here, we want to give some concrete examples of functional spaces for which the theory holds. Beyond the classical $L^p$-spaces, other remarkable functional spaces arising from suitable convex $\varphi$-functions can be included. For example, \begin{equation}\label{zygmund}\varphi_{\beta,\gamma}(u):=u^\beta\log^\gamma(u+e),\quad u\ge 0,\end{equation} where $\beta\ge 1$, $\gamma>0$, and \begin{equation}\label{exp}\varphi_\rho(u):=e^{u^{\rho}}-1,\quad u\ge 0,\end{equation} with $\rho>0$, generate the $L^\beta \log^\gamma L$-spaces (or \textit{Zygmund spaces}) and the \textit{exponential spaces}, respectively. Particularly, Zygmund spaces are extensively employed in the theory of partial differential equations, while those of exponential-type are used for establishing embedding theorems (see, e.g., \cite{stein69,BR80} and \cite{H03}, respectively). The convex modular functionals associated with $\varphi_{\beta,\gamma}$ and $\varphi_\rho$ are represented by
\[
I^{\varphi_{\beta,\gamma}}[f]:=\int_I\left|f(x)\right|^\beta\log^\gamma(e+\left|f(x)\right|)dx,\quad f\in M(I),
\]
and
\[
I^{\varphi_\rho}[f]:=\int_I\left(e^{|f(x)|^\rho}-1\right)dx,\quad f\in M(I),
\]
respectively.
\vskip0.2cm

  Since we have to assume that the above $\varphi$-functions are also $N$-functions, in the case of Lebesgue spaces we must consider $\varphi(u)=u^p$ with $1<p<+\infty$, in the case of Zygmund spaces we take $\varphi_{\beta,\gamma}$ given in (\ref{zygmund}) with $\beta >1$, as well as for exponential spaces we consider $\varphi_\rho$ introduced in (\ref{exp}) with $\rho\ge2$.
\vskip0.2cm

  In literature, the class of $N$-functions holds a central position. Indeed, the latter are used to generate special instances of Orlicz spaces that are useful for several applications, especially in statistical mechanics. In \cite{IE1972}, the following strongly nonlinear $N$-function was used
\begin{equation}\label{exp2}
\varphi(u):=e^{u}-u-1,\quad u\ge 0.
\end{equation}
More recently, in \cite{ML2014}, this approach was adopted to describe the statistics of large regular statistical systems, both classical and quantum, using a pair of Orlicz spaces, namely $L^{\cosh-1}$ and $L\log(L+1)$.
The first Orlicz space $L^{\cosh-1}$, generated by the $N$-function \begin{equation}\label{exp3}\varphi(u):=\cosh u-1,\quad u\ge 0,\end{equation} appears as a suitable framework for describing the set of regular observable. The second Orlicz space $L\log(L+1)$ is a Zygmund space defined by the $N$-function $u\log(u+1)$, $u\ge0$. This approach is particularly suitable for large systems, i.e., systems with infinite degrees of freedom that demand more general Banach spaces than $L^p$-spaces.
\vskip0.2cm

  Now, we point out that the theory presented here can be applied to several sigmoidal functions $\sigma$ (even not necessarily belonging to $C^2(\R)$) generating compactly supported density functions (for more details, see \cite{coroianu2024approximation}). For instance, the well-known \textit{ramp function} $\sigma_r$ is defined by
\[
\sigma_r(x):=\begin{cases} 0, &\mbox{ $x<-3/2$} \\ x/3+1/2, &\mbox{ $-3/2<x<3/2$}\\ 1, & \mbox{ $x>3/2$,}
\end{cases}
\]
 and the corresponding $\phi_{\sigma_r}$ is a compactly supported function.
\vskip0.2cm

  Other sigmoidal functions can be generated by the well-known \textit{central B-splines of order $n\in\mathbb{N}$}, given by
\begin{equation*}\label{spline}
\beta_n(x):=\frac{1}{(n-1)!}\sum_{j=0}^n(-1)^j\binom{n}{j}\left(\frac{n}{2}+x-j\right)_+^{n-1},\qquad x\in\mathbb{R},
\end{equation*}
where $(\cdot)_+$ denotes the positive part, i.e., $(x)_+:=\max\left\{x,0\right\}$, that lead to
\[
\sigma_{\beta_i}(x):=\int_{-\infty}^x\beta_s(t)dt,\qquad x\in\mathbb{R},
\]
whose support is contained in $[k,+\infty)$ (with $k>0$ is a suitable constant).
\vskip0.2cm

  The presented theory is not limited to the sigmoidal activation functions but also includes the widely used ReLU activation function (Rectified Linear Unit) (see, e.g., \cite{DDFH21,OSZ21,Y21,plonka2023spline,LMP24}). Additionally, we can consider its $k$-th power, referred to as ReLU$^{k}$, which is also known as Rectified Power Units (RePUs) (see \cite{LTY19} and \cite{GR21}). Interested readers can find further details in \cite{CJat23}.

Now, we give a quantitative result which follows from Theorem \ref{newthm}.

\begin{theorem}\label{mainthm2orlicz}
Let $\sigma$ be a sigmoidal function and $\varphi$ be a convex $\varphi$-function satisfying the $\Delta^\prime$-condition and such that $u\mapsto u^{-\beta}\varphi(u)$, $u>0$, is increasing, with $\beta>1$. Moreover, let $M_{0,1}^\varphi(\phis)<+\infty$. Thus, for every $f\in L^{\varphi}(I)$, it turns out that
\begin{equation*}
\| K_nf-f\|_\varphi\le\Psi^{\varphi,\sigma} \cdot \omega\left( f,\frac{1}{n}\right)_\varphi,
\end{equation*}
where $\Psi^{\varphi,\sigma}>0$ is a suitable constant.
\end{theorem}
\noindent The proof is the same of Theorem \ref{mainthm1orlicz}, where we use Theorem \ref{newthm} in place of Theorem \ref{thm1orlicz}.
\vskip0.2cm

  We point out that both the remarkable cases of $L^p$-spaces and Zygmund spaces $L^\beta\log^\gamma$ are included in Theorem \ref{thm2}, since the corresponding $\varphi$-functions satisfy $\Delta^\prime$-condition (see, e.g., \cite{krivoshein2022wavelet}). Moreover, the moment condition $M_{0,1}^\varphi(\phis)<+\infty$ is satisfied if the decay condition (\ref{a2}) holds for $\alpha>p$ in case of $\varphi(u)=u^p$ or $\alpha>\beta+\gamma$ in case of $\varphi_{\beta,\gamma}(u):=u^\beta\log^\gamma(u+e)$. This is true if we require that condition ($S3$) holds for $\alpha>p$ or $\alpha>\beta+\gamma$, respectively. 
\vskip0.2cm 

  As an example of sigmoidal function, we may recall a particular instance depending on a parameter $\theta>0$, defined as
\[
\sigma_\theta(x):=\displaystyle \begin{cases} \displaystyle\frac{2^{-\theta}}{4}|x|^{-\theta}, &\mbox{ $x\le-\frac{1}{2}$} \\ \displaystyle\frac{1}{2}x+\frac{1}{2}, &\mbox{ $-\frac{1}{2}<x\le\frac{1}{2}$}\\ 1-\displaystyle\frac{2^{-\theta}}{4}|x|^{-\theta}, & \mbox{ $x>\frac{1}{2}$,}
\end{cases}
\]
with the corresponding density function given by
\[
\phi_{\sigma_\theta}(x)\, :=\, \frac{1}{2}[\sigma_\theta(x+1)-\sigma_\theta(x-1)], \hskip1cm x \in \R,
\]
(see, e.g., \cite{CC24}). The above sigmoidal function allows us to calibrate the parameter $\theta$ based on the considered functional spaces: in particular, using the $L^p$-spaces, we should require that $\theta>p+1$, while for Zygmund spaces we need to assume $\theta>\beta+\gamma+1$.
\vskip0.2cm

  As a consequence of Theorem \ref{mainthm1orlicz} and Theorem \ref{mainthm2orlicz}, we are able to state a strong convergence theorem, i.e, a Luxemburg norm convergence theorem, which extends Theorem 4.1 of \cite{CVMathNachr} from the space of bounded and uniformly continuous functions on $I$ to the whole Orlicz space $L^\varphi(I)$ (so as to include also not necessarily continuous functions).
\begin{theorem}\label{luxemburgconvergence}
Let $\sigma$ be a sigmoidal function and $f\in L^{\varphi}(I)$.
\begin{enumerate}
    \item If $\varphi$ is an $N$-function and $\phi_\sigma$ has compact support, then
\begin{equation*}
\lim_{n\to+\infty}\| K_nf-f\|_\varphi= 0.
\end{equation*}

\item If $\varphi$ is a convex $\varphi$-function satisfying the $\Delta^\prime$-condition, such that $u\mapsto u^{-\beta}\varphi(u)$, $u>0$, is increasing, with $\beta>1$, then
\begin{equation*}
\lim_{n\to+\infty}\| K_nf-f\|_\varphi= 0.
\end{equation*}
\end{enumerate}
\end{theorem}

The convergence result stated above is stronger than the modular convergence theorem provided in Theorem 4.4 of \cite{CVMathNachr}. Indeed, it is a well-known fact that Luxemburg norm convergence implies modular convergence, as extensively discussed in Section \ref{sec3}. However, if in Theorem 4.4 of \cite{CVMathNachr} we require (the additional assumption) that $\varphi$ satisfies the $\Delta_2$-condition, it can be also interpreted as a Luxemburg norm convergence theorem. But in this case, the required additional assumption excludes the application to certain noteworthy examples of Orlicz spaces, such as the exponential spaces. Now, by Theorem \ref{luxemburgconvergence} we get the Luxemburg norm convergence of the Kantorovich NN operators without requiring the $\Delta_2$-condition, and then it applies also, e.g., for the exponential-type spaces (which are generated by N-functions).
\vskip0.2cm

In order to establish a weaker version of the quantitative estimates achieved in Theorem \ref{mainthm1orlicz} and in Theorem \ref{mainthm2orlicz}, we need to use a modular version of the involved modulus of smoothness. Therefore, we introduce the \textit{weak $\varphi$-modulus of smoothness} in Orlicz spaces, namely
\begin{equation}\label{wmodulusofsmoothnss}
\widetilde{\omega}(f,\delta)_\varphi:=\sup_{|h|\le\delta} I^\varphi\left[f(\cdot+h)-f(\cdot) \right],\qquad\delta>0.
\end{equation}
\noindent Also here, to avoid possible issues, we can extend the involved function $f$ with periodicity on the whole $\R$. Furthermore, we recall the following fundamental fact. It is well-known that, for any $f\in L^\varphi(I)$, there exists a suitable $\lambda^\star>0$ such that 
\begin{equation}\label{property}
\widetilde{\omega}(\lambda^\star f,\delta)_\varphi\to 0,\quad \delta\to 0^+,
\end{equation}
(see Theorem 2.4 of \cite{BMV2003}).
\vskip0.2cm

By using (\ref{wmodulusofsmoothnss}), we can give the following weak version of Theorem \ref{SPOrlicz}.

\begin{lemma}\label{lemmaweak} Let $f\in L^\varphi(I)$ with $\varphi$ a convex $\varphi$-function. Thus, there exists $f_{1,h}\in W^{1,\varphi}(I)$, i.e., the first-order Steklov function given in (\ref{steklov}), such that the following holds:
\begin{itemize}
    \item [(i)] $f_{1,h}\in\mathcal{W}^{1,\varphi}$, for every fixed $0<h\le \bar{h}$, with $\bar{h}$ sufficiently small;
    \item [(ii)] there exists $\lambda>0$ such that 
    \[
I^\varphi\left[\lambda\left(f-f_{1,h}\right)\right]\le \widetilde{\omega}(\lambda f,h)_\varphi, \quad 0<h\le b-a;
\]
 \item [(iii)] for any $\lambda>0$, there holds \[I^\varphi\left[\lambda f^\prime_{1,h}\right]\le \widetilde{\omega}\left(\frac{\lambda}{h} f,h\right)_\varphi,\quad 0<h\le b-a;\]
  \item [(iv)] for any $\lambda>0$, there holds \[I^\varphi\left[\lambda \left(f^\prime_{1,h}\right)^2\right]\le \widetilde{\omega}\left(\frac{\lambda}{h^2} f,h\right)_\varphi,\quad 0<h\le \bar{h}.\]
\end{itemize}
 \end{lemma}
\begin{proof}
Let $x\in I$, $n\in\N$ and $$f_{1,h}(x)=\displaystyle\frac{1}{h}\int_0^h f(x+t) dt,$$ be the first-order Steklov function given in (\ref{steklov}). Similarly as in (ii) of Theorem \ref{SPOrlicz}, one can see that $f_{1,h}\in W^{1,\varphi}(I)$. In particular, by a suitable change of variable, we can note that $f_{1,h}\in AC(I)$.
\vskip0.2cm

  Now, we want to prove (i). Since $f\in L^\varphi(I)$, there is $\lambda^\star>0$ such that (\ref{property}) holds. Firstly, we need to observe that being \(f \in L^\varphi(I)\), it turns out that
\begin{equation}\label{absurd}
\left|f(x+h)-f(x)\right|\le 1
\end{equation}
for \(x \in I\) a.e. and $h>0$ sufficiently small. We will prove this by contradiction. Denoting by $\mu$ the Lebesgue measure, we suppose there exists a subset \(\Omega \subset I\), with \(\mu(\Omega) > 0\), such that
\[
\left|f(x+h)-f(x)\right| > 1
\]
for \(x \in \Omega\). Consequently, we have
\[
\lambda^\star \left|f(x+h)-f(x)\right| > \lambda^\star
\]
for \(x \in \Omega\), where \(\lambda^\star>0\) is defined above. Considering the non-decreasing property of the $\varphi$-function, we can write
\[
I^\varphi\left[\lambda^\star\left(f(\cdot+h)-f(\cdot)\right)\right]\geq\int_\Omega \varphi\left(\lambda^\star \left|f(x+h) - f(x)\right|\right) \, dx > \int_\Omega \varphi\left(\lambda^\star\right) \, dx = \mu(\Omega) \cdot \varphi(\lambda^\star) > 0,
\]
from which we get a contradiction as \(h \to 0^+\).
Now, we recall that
\[
f_{1,h}^\prime(x)=\frac{1}{h}\left(f(x+h)-f(x)\right),\qquad h>0,\ x\in I.
\]
One can see that  $f_{1,h}\in \mathcal{W}^{1,\varphi}(I)$, for $h>0$ sufficiently small. Indeed, by using (\ref{absurd}) and the absolute homogeneity of $\|\cdot\|_\varphi$, we have that 
\begin{equation*}
\begin{split}
\|\left(f_{1,h}^\prime\right)^2\|_\varphi&=\left\|\frac{1}{h^2}\left(f(\cdot+h)-f(\cdot)\right)^2\right\|_\varphi
\\
&=\frac{1}{h^2}\cdot\inf\left\{u>0 \ : \ I^\varphi\left[\frac{\left(f(\cdot+h)-f(\cdot)\right)^2}{u}\right]\le 1\right\}
\\
&\le \frac{1}{h^2}\cdot\inf\left\{u>0 \ : \ I^\varphi\left[\frac{f(\cdot+h)-f(\cdot)}{u}\right]\le 1\right\}
\\
&=\frac{1}{h^2}\cdot \|f(\cdot+h)-f(\cdot)\|_\varphi
\\
&\le \frac{1}{h^2}\cdot \omega\left(f,h\right)_\varphi<+\infty,
\end{split}
\end{equation*}
for every fixed $h>0$ sufficiently small, since $f\in L^\varphi(I)$.
\\
Now, we prove (ii). To this aim, we consider
\begin{equation*}
\begin{split}
f(x)-f_{1,h}(x)&=f(x)-\frac{1}{h}\int_0^h f(x+t) \ dt
\\
&=\frac{1}{h}\int_0^h\left\{ f(x)-f(x+t)\right\} \ dt.
\end{split}
\end{equation*}
Now, we fix $\lambda>0$, such that $\lambda<\lambda^\star$, where $\lambda^\star$ is defined above. Passing to the modular and using the weak Minkowski inequality (\ref{weakminkineq}) with $g\equiv 1$, it turns out that
\begin{equation*}
\begin{split}
I^\varphi\left[\lambda\left(f-f_{1,h}\right)\right]&= \int_a^b\varphi\left(\lambda\left|f(x)-f_{1,h}(x)\right|\right)\ dx
\\
&=\int_a^b\varphi\left(\frac{\lambda}{h}\int_0^h\left|f(x)-f(x+t)\right|dt\right) dx
\\
&\le\frac{1}{h}\int_0^h\left[\int_a^b\varphi\left( \lambda\left|f(x)-f(x+t)\right|\right)dx\right] dt
\\
&\le \widetilde{\omega}(\lambda f, h)_\varphi
\le \widetilde{\omega}(\lambda^\star f, h)_\varphi<+\infty,
\end{split}
\end{equation*}
in view of Theorem 2.2 (c) of \cite{BMV2003}.
\vskip0.2cm

  Now, we prove (iii). Considering $\lambda>0$ and $0<h\le b-a$, we immediately see that
\begin{equation*}
\begin{split}
I^\varphi\left[\lambda f_{1,h}^\prime\right]&= \int_a^b\varphi\left[\frac{\lambda}{h}\left(f(x+h)-f(x)\right)\right]dx
\le \widetilde{\omega}\left(\frac{\lambda}{h} f, h\right)_\varphi.
\end{split}
\end{equation*}
Finally, to prove (iv) we proceed similarly as in (iii). Taking $\lambda>0$ and $0<h\le b-a$, by using again (\ref{absurd}), we get
\begin{equation*}
\begin{split}
I^\varphi\left[\lambda \left(f_{1,h}^\prime\right)^2\right]&= \int_a^b\varphi\left[\frac{\lambda}{h^2}\left(f(x+h)-f(x)\right)^2\right]dx
\\
&\le \int_a^b\varphi\left[\frac{\lambda}{h^2}\left(f(x+h)-f(x)\right)\right]dx
\le \widetilde{\omega}\left(\frac{\lambda}{h^2} f, h\right)_\varphi.
\end{split}
\end{equation*}
\end{proof}
\begin{remark}\rm
From (i) and (ii) of Lemma \ref{lemmaweak}, it turns out that the space $\mathcal{W}^{1,\varphi}(I)$ is modularly dense in $L^\varphi(I)$.
\end{remark}
\begin{remark} \rm
 We point out that in (iii) and (iv) of Lemma \ref{lemmaweak}, the right-hand side of the inequalities is not necessarily finite. To have this, in the next theorems additional assumptions are needed.
\end{remark}

Now, we want to obtain the quantitative version of Theorem \ref{thm3}, based on the weak $\varphi$-modulus of smoothness defined in (\ref{wmodulusofsmoothnss}).
\begin{theorem}\label{mainthm3orlicz} Let $\sigma$ be a sigmoidal function and $\varphi$ be a convex $\varphi$-function such that $u\mapsto u^{-\beta}\varphi(u)$, $u>0$, is increasing, with some $\beta>1$. Moreover, let $M_{0,2}^\varphi(\phis)<+\infty$. Thus, for every $f\in L^{\varphi}(I)$, there exists $\lambda>0$ such that for every sufficiently large $n\in\N$ there holds
\begin{equation*}\label{est1}
I^\varphi\left[\lambda\left( K_nf-f\right)\right]\le\Xi^{\sigma,\varphi} \cdot \widetilde{\omega}\left(12\lambda L^2 f,\frac{1}{\sqrt{n}}\right)_\varphi+\xi_1^{\varphi,\sigma}\cdot\frac{3\lambda}{n},
\end{equation*}
where $\Xi^{\sigma,\varphi}>0$ is a suitable constant, $L\ge1$, as stated in Remark \ref{rmkhl}, arises from (\ref{equivalence}) and $\xi_1^{\sigma,\varphi}>0$ arises from Theorem \ref{thm3}.
\end{theorem}
\begin{proof}
Let $f\in L^{\varphi}(I)$ be fixed (and extended on $\R$ as a $(b-a)$-periodic function). Since $f\in L^\varphi(I)$, there is $\lambda^\star>0$ such that (\ref{property}) holds. Therefore, let $\lambda>0$ be such that $12\lambda L^2<\lambda^\star$.\\ By Lemma \ref{lemmaweak} (ii), there exists $f_{1,h}\in W^{1,\varphi}(I)$, $h>0$, such that
\[
I^\varphi\left[\lambda\left(f-f_{1,h}\right)\right]\le \widetilde{\omega}(\lambda f,\delta)_\varphi.
\]
Moreover, by Lemma \ref{lemmaorlicz}, we know that
\begin{equation}\label{est2}
I^\varphi\left[ \lambda\left(K_n\left(f-f_{1,h}\right)\right)\right]\le \phis(2)^{-1} I^\varphi\left[\lambda\left(f-f_{1,h}\right)\right],
\end{equation}
being $f-f_{1,h}\in L^\varphi(I)$.\\
Then, by the convexity of the modular, the linearity of NN Kantorivich operators and (\ref{est2}), we obtain
\begin{equation*}
\begin{split}
I^\varphi\left[\lambda\left( K_nf-f\right)\right]&\le I^\varphi\left[3\lambda\left( K_nf-K_nf_{1,h}\right)\right]+I^\varphi\left[3\lambda\left(K_nf_{1,h}-f_{1,h}\right)\right]+I^\varphi\left[3\lambda\left(f_{1,h}-f\right)\right]
\\
&\le \left(\phis(2)^{-1}+1\right) I^\varphi\left[3\lambda \left(f_{1,h}-f\right)\right]+I^\varphi\left[3\lambda\left(K_nf_{1,h}-f_{1,h}\right)\right]
\\
&\le \left(\phis(2)^{-1}+1\right)\widetilde{\omega}(3\lambda f,h)_\varphi +I^\varphi\left[3\lambda\left(K_nf_{1,h}-f_{1,h}\right)\right].
\end{split}
\end{equation*}

Hence, by Theorem \ref{thm3}, (i) and (iv) of Lemma \ref{lemmaweak}, taken $h=\frac{1}{\sqrt{n}}\le b-a$, for $n\in\N$ sufficiently large, we achieve
\begin{align}\label{min1}
I^\varphi\left[3\lambda\left(K_nf_{1,\frac{1}{\sqrt{n}}}-f_{1,\frac{1}{\sqrt{n}}}\right)\right]&\le \xi_1^{\sigma,\varphi}\cdot\frac{3\lambda}{n}+\xi_2^{\sigma,\varphi}\cdot I^\varphi\left[\frac{12\lambda L^2}{n} \left(f_{1,\frac{1}{\sqrt{n}}}^\prime\right)^2\right]\nonumber
\\
&\le   \xi_1^{\sigma,\varphi}\cdot\frac{3\lambda}{n}+\xi_2^{\sigma,\varphi}\cdot\widetilde{\omega}\left(\frac{12\lambda L^2}{n}\cdot n f,\frac{1}{\sqrt{n}}\right)_\varphi\nonumber
\\
&\le  \xi_1^{\sigma,\varphi}\cdot\frac{3\lambda}{n}+\xi_2^{\sigma,\varphi}\cdot\widetilde{\omega}\left(\lambda^\star f,\frac{1}{\sqrt{n}}\right)_\varphi<+\infty,\nonumber
\end{align}
where $\xi_1^{\sigma,\varphi}$ and $\xi_2^{\sigma,\varphi}$ are those of Theorem \ref{thm3}.
\\
Finally, being $L\ge1$ and using Theorem 2.2 (c) of \cite{BMV2003}, we have
\begin{equation*}
\begin{split}
I^\varphi\left[\lambda\left( K_nf-f\right)\right]&\le \left(\phis^{-1}(2)+1\right)\widetilde{\omega}\left(3\lambda f,\frac{1}{\sqrt{n}}\right)_\varphi+\xi_1^{\sigma,\varphi}\cdot\frac{3\lambda}{n}+\xi_2^{\sigma,\varphi}\cdot  \widetilde{\omega}\left(12\lambda L^2 f,\frac{1}{\sqrt{n}}\right)_\varphi
\\
&\leq\left(\phis^{-1}(2)+1+\xi_2^{\sigma,\varphi}\right)\cdot\widetilde{\omega}\left(12\lambda L^2 f,\frac{1}{\sqrt{n}}\right)_\varphi+\xi_1^{\sigma,\varphi}\cdot\frac{3\lambda}{n}
\\
&=:\Xi^{\sigma,\varphi}\cdot\widetilde{\omega}\left(12\lambda L^2 f,\frac{1}{\sqrt{n}}\right)_\varphi+\xi_1^{\sigma,\varphi}\cdot\frac{3\lambda}{n}.
\end{split}
\end{equation*}
This concludes the proof.
\end{proof}

  The previous estimate holds in every type of Sobolev-Orlicz spaces, included those generated by exponential $\varphi$-function such as (\ref{exp}), (\ref{exp2}) or (\ref{exp3}). In order to have the finiteness of the discrete absolute $\varphi$-moment $M_{0,2}^\varphi(\phis)$, we need to use sigmoidal functions with exponential decay as the variable approaches $-\infty$. To this regard, we may recall the well-known \textit{logistic function} defined by
\begin{equation}\label{logisticf}
\sigma_l(x):=(1+e^{-x})^{-1},\qquad x\in\mathbb{R},
\end{equation}
and the \textit{hyperbolic tangent function}, given by
\[
\sigma_h(x):=(\tanh x+1)/2,\qquad x\in\mathbb{R}.
\]
Both of them easily satisfy condition ($S 3$) for every $\alpha>0$, in view of their exponential decay to zero, as $x\to-\infty$, together with ($S_1$) and ($S_2$).
\vskip0.2cm

  In particular, taken as $\varphi$-function that given in (\ref{exp}), we may consider the logistic function (\ref{logisticf}) as the sigmoidal activation function. Therefore, by choosing $\rho<1/2$, where $\rho>0$ is the parameter of the exponential $\varphi$-function given in (\ref{exp}), one can easily see that the moment assumption $M_{0,2}^\varphi(\phis)<+\infty$ turns out to be certainly satisfied.
\begin{remark}
    \rm We highlight that Theorem \ref{mainthm3orlicz} can be read as a modular convergence theorem on the whole $L^\varphi(I)$, where $\varphi$ is simply a convex $\varphi$-function.
\end{remark}
Furthermore, as a consequence of Theorem \ref{weakthm1orlicz}, we may give another weak quantitative estimate.
\begin{theorem}\label{secondmodestimate} Let $\varphi$ be a convex $\varphi$-function and $\sigma$ be a sigmoidal function such that the resulting $\phi_\sigma$ has compact support contained in $[-\Upsilon, \Upsilon]\subset\mathbb{R}$, with $\Upsilon>0$. Thus, for every $f\in L^\varphi(I)$, there exists $\lambda>0$ such that
\[
I^\varphi[\lambda\left(K_nf-f\right)]\le (2+\phis(2)^{-1})\cdot\widetilde{\omega}\left(6\lambda(1+\Upsilon)f, \frac{1}{n}\right)_\varphi,
\]
for every sufficiently large $n\in\N$.
\end{theorem}
\begin{proof}
  Let $f\in L^{\varphi}(I)$ be fixed (and extended on $\R$ as a $(b-a)$-periodic function). Since $f\in L^\varphi(I)$, there is $\lambda^\star>0$ such that (\ref{property}) holds. Therefore, let $\lambda>0$ be such that $6\lambda (1+\Upsilon)<\phis(2)\lambda^\star$.\\ Proceeding as in the first part of the proof of Theorem \ref{mainthm3orlicz}, by Lemma \ref{lemmaorlicz}, Theorem \ref{weakthm1orlicz} and Lemma \ref{lemmaweak} (ii), there exists $f_{1,h}\in W^{1,\varphi}(I)$, $h>0$, such that
  \begin{equation*}
      \begin{split}
          I^\varphi[\lambda\left(K_nf-f\right)]&\le(1+\phis(2)^{-1})\cdot\widetilde{\omega}(3\lambda f, h)_\varphi+I^\varphi\left[3\lambda\left(K_nf_{1,h}-f_{1,h}\right)\right]
          \\
          &\le (1+\phis(2)^{-1})\cdot\widetilde{\omega}(3\lambda f, h)_\varphi+I^\varphi\left[\frac{6\lambda(1+\Upsilon)}{n}f_{1,h}^\prime\right].
      \end{split}
  \end{equation*}
  Finally, taking $h=\frac{1}{n}\le b-a$, for $n\in\N$ sufficiently large, and using Lemma \ref{lemmaweak} (iii), we achieve the thesis, i.e.,
  \begin{equation*}
      \begin{split}
          I^\varphi[\lambda\left(K_nf-f\right)]&\le(1+\phis(2)^{-1})\cdot\widetilde{\omega}\left(3\lambda f, \frac{1}{n}\right)_\varphi+\widetilde{\omega}\left(6\lambda(1+\Upsilon) f, \frac{1}{n}\right)_\varphi
          \\
          &\le (2+\phis(2)^{-1})\cdot\widetilde{\omega}\left(6\lambda(1+\Upsilon)f, \frac{1}{n}\right)_\varphi.
      \end{split}
  \end{equation*}
\end{proof}
\begin{remark} \rm
    Under the assumptions of Theorem \ref{secondmodestimate}, a further modular convergence theorem can be derived for the entire $L^\varphi(I)$, where the involved $\varphi$-function can now be taken without restrictions. Therefore, the delicate $L^1$-case, i.e., when $\varphi(u) = u$, $u \ge 0$, can also be included. Since such a function satisfies the $\Delta_2$-condition, the convergence is actually strong. Moreover, we highlight that Theorem \ref{secondmodestimate} improves Theorem 4.7 of \cite{cpNNlp} (when $p=1$) in terms of the rate of convergence. This improvement is essentially due to the current assumption that the density function $\phis$ is compactly supported.
        \end{remark}

\section{Qualitative analysis in Orlicz Lipschitz classes}\label{s6}

From the above quantitative estimates we can directly deduce the qualitative rate of convergence, assuming that $f$ belongs to suitable Lipschitz classes defined in the Orlicz-setting. To this aim, we recall the definition of \textit{Lipschitz classes} in terms of the $\varphi$-modulus of smoothness. Such classes are defined as
\begin{equation}\label{Lipphi}
Lip(\nu,\varphi):=\left\{ f\in L^\varphi(I):\omega(f,\delta)_\varphi=\mathcal{O}(\delta^\nu),\,\text{ as }\delta\rightarrow 0^+\right\},
\end{equation}
with $0<\nu\le1$.
As a consequence of Theorem \ref{mainthm1orlicz} and Theorem \ref{mainthm2orlicz}, we obtain the following.
\begin{corollary}\label{cordirectresult} Let $f\in Lip(\nu,\varphi)$, with $0<\nu\le1$. Under the assumptions made in Theorem \ref{mainthm1orlicz} (or in Theorem \ref{mainthm2orlicz}), there exists $\Gamma>0$ such that
\[
\| K_nf-f\|_\varphi\le \Gamma n^{-\nu},
\]
for sufficiently large $n\in\N$, with $\Gamma:=\Lambda C$ where $\Lambda$ and $C$ are suitable constants arising from Theorem \ref{mainthm1orlicz} (or Theorem \ref{mainthm2orlicz}, respectively) and (\ref{Lipphi}), respectively. 
\end{corollary}

Finally, we may also define weak Lipschitz classes, where in (\ref{Lipphi}) the $\varphi$-modulus of smoothness is replaced by its weak version, namely
\begin{equation}\label{Lipphiweak}
\widetilde{Lip}(\nu,\varphi):=\left\{ f\in L^\varphi(I): \exists \ \lambda>0 \ \text{ s.t. } \ \widetilde{\omega}(\lambda f,\delta)_\varphi=\mathcal{O}(\delta^\nu),\,\text{ as }\delta\rightarrow 0^+\right\},
\end{equation}
with $0<\nu\le1$. \\
Now, it seems natural to ask something about the relation between (\ref{Lipphi}) and (\ref{Lipphiweak}). Hence, we can prove the following.
\begin{proposition}\label{inclusion} Let $\varphi$ be a convex $\varphi$-function and $0<\nu\le 1$. Thus, there holds
\[
Lip(\nu,\varphi)\subseteq \widetilde{Lip}(\nu,\varphi).
\]
\end{proposition}
\begin{proof}
Let $f\in Lip(\nu,\varphi)$, with $0<\nu\le 1$. By definition, there exist a constant $C>0$ and $\bar{\delta}>0$ such that $\omega(f,\delta)_\varphi\le C\delta^\nu<1$, for $0<\delta\le \bar{\delta}$. In particular, we also have
\[
\|f(\cdot+t)-f(\cdot)\|_\varphi\leq C \delta^\nu<1,
\]
for every $t\in\R$ such that $|t|<\delta$.
By using Theorem 1.1 (c) of \cite{BMV2003}, we obtain
\[
I^\varphi\left[f(\cdot+t)-f(\cdot)\right]\le \|f(\cdot+t)-f(\cdot)\|_\varphi\leq C\delta^\nu,
\]
for every $t\in\R$ such that $|t|<\delta$. Passing now to the supremum with respect to $t$, we get $\widetilde{\omega}(\lambda f,\delta)_\varphi\leq C\delta^\nu$, for $0<\lambda\le 1$, and $0<\delta<\bar{\delta}$. This concludes the proof.
\end{proof}
Now, for completeness, it is important to provide instances of $\varphi$-functions for which the inclusion shown in Proposition \ref{inclusion} is proper. We will illustrate this with an example, but first we need some general theoretical results.
\begin{proposition} Let $f\in L^\varphi(\Omega)$, $\Omega\subseteq\R$, with $\varphi$ satisfying the $\Delta_2$-condition, and $0<\nu\le 1$. If there exists $\lambda^\star>0$ such that 
\begin{equation*}
    \widetilde{\omega}\left(\lambda^\star f,\delta\right)=\mathcal{O}(\delta^\nu),\quad \delta\to 0^+,
\end{equation*}
    then
    \begin{equation*}
    \widetilde{\omega}\left(\lambda f,\delta\right)=\mathcal{O}(\delta^\nu),\quad \delta\to 0^+,
\end{equation*}    
for every $\lambda>0$. Obviously, also the converse is always trivially true.
\end{proposition}
\begin{proof}
    Let $f\in L^\varphi(\Omega)$ and $\lambda>0$ be fixed. If $\lambda\leq\lambda^\star$, then 
    \[
    \widetilde{\omega}\left(\lambda f,\delta\right)\le \widetilde{\omega}\left(\lambda^\star f,\delta\right)=\mathcal{O}(\delta^\nu),\quad \delta\to 0^+,
    \]
    in view of property (a) of Theorem 2.2 of \cite{BMV2003}. Otherwise, if $\lambda>\lambda^\star$, then there exists $N\in\N$ such that $\lambda\leq 2^N\lambda^\star$, from which
    \[
    \widetilde{\omega}\left(\lambda f,\delta\right)\le \widetilde{\omega}\left(2^N\lambda^\star f,\delta\right)\leq C_\varphi^N\widetilde{\omega}\left(\lambda^\star f,\delta\right)=\mathcal{O}(\delta^\nu),\quad \delta\to 0^+,
    \]
    where we use, by easy induction, that $\varphi(2^Nu)\leq C_\varphi^N\varphi(u)$, for $u\geq 0$ and $C_\varphi>0$ by the $\Delta_2$-condition in (\ref{delta2}).
\end{proof}
\begin{proposition}\label{noLip}
    Let $f\in L^\varphi(\Omega)$, $\Omega\subseteq\R$. If there exists $\lambda^\star>0$ such that $\widetilde{\omega}(\lambda^\star f,\delta)\nrightarrow 0$ as $\delta\to 0^+$, then
    \[
    f\notin Lip(\nu,\varphi),\qquad 0<\nu\le 1.
    \]
\end{proposition}
\begin{proof}
    Let $\lambda^\star>0$ be such that $\widetilde{\omega}(\lambda^\star f,\delta)\nrightarrow 0$ as $\delta\to 0^+$. Assuming by contradiction that $f\in Lip(\nu,\varphi)$, with $0<\nu\leq1$, there exist $C>0$ and $\bar{\delta}_{\lambda^\star}>0$ such that $\lambda^\star\omega(f,\delta)\leq C\lambda^\star\delta^\nu<1$, for every $0<\delta\leq\bar{\delta}_{\lambda^\star}$, meaning that
    \[
    \lambda^\star\left\|f(\cdot+t)-f(\cdot)\right\|_\varphi\leq \lambda^\star C \delta^\nu<1,
    \]
    for every $|t|<\delta$, $0<\delta\leq \bar{\delta}_{\lambda^\star}$. By using again Theorem 1.1 (c) of \cite{BMV2003}, we obtain
    \[
    I^\varphi\left[\lambda^\star\left(f(\cdot+t)-f(\cdot)\right)\right]\leq\lambda^\star\left\|f(\cdot+t)-f(\cdot)\right\|_\varphi\leq \lambda^\star C\delta^\nu,
    \]
    for every $|t|<\delta$, $0<\delta\leq \bar{\delta}_{\lambda^\star}$. Passing to the supremum with respect to $t\in\R$ such that $|t|<\delta$, we get
    \[
    0\leq\widetilde{\omega}(\lambda^\star f,\delta)\leq \lambda^\star\omega(f,\delta)\leq \lambda^\star C\delta^\nu.
    \]
    Hence, considering the limit as $\delta\to 0^+$ in the above sequence of inequalities, we obtain that $\widetilde{\omega}(\lambda^\star f,\delta)\to 0$, which is a contradiction.
\end{proof}
Now, we are ready to give a concrete example showing that the inclusion stated in Proposition \ref{inclusion} is, in general, proper.
\begin{example}\rm Let $\Omega:=[0,1]$. We consider $f(x):=\ln \left(x^{-1/2}\right)$, if $x\in(0,1]$ and $f(0):=0$, and we extend it on the whole $\R$ as a $1$-periodic function. We take into account the exponential $\varphi$-function $\varphi(u)=e^u-1$, $u\geq 0$, given in (\ref{exp}), which clearly does not satisfy the $\Delta_2$-condition.
\\
Now, for every fixed $t\in [0,1/2]$, we want to estimate the modular $I^\varphi\left[\lambda\left(f(\cdot+t)-f(\cdot)\right)\right]$, considering different values of $\lambda>0$. 
\\
Firstly, we take $\lambda=1$. Thus, by the $1$-periodicity of $f$, we can write
\begin{equation*}
    \begin{split}
        \int_\Omega\varphi\left(\left|f(x+t)-f(x)\right|\right) dx&=\int_0^{1-t}\varphi\left(\left|f(x+t)-f(x)\right|\right) dx+\int_{1-t}^1\varphi\left(\left|f(x-1+t)-f(x)\right|\right) dx
	\\&=: T_1(t)+T_2(t),
 \end{split}
\end{equation*}
where $t\in\left[0,1/2\right]$.
\\
We now focus on $T_1$. By the change of variable $z=x+t$, we get
\begin{equation*}
    \begin{split}
        T_1(t)&= \int_0^{1-t}\varphi\left(\left|f(x+t)-f(x)\right|\right) dx
        \\
        &=\int_t^1\varphi\left(\left|f(z)-f(z-t)\right|\right) dz 
        \\
        &=\int_t^1\left(\sqrt{\frac{z}{z-t}}-1\right)dz 
        \end{split}
\end{equation*}
\begin{equation*}
    \begin{split}
        &=\lim_{y\to t}\left[z\sqrt{1 - \frac{t}{z}}  - \frac{1}{2} t \log\left(1 - \sqrt{1 - \frac{t}{z}}\right) + \frac{1}{2} t \log\left(1 + \sqrt{1 - \frac{t}{z}}\right)-z\right]_y^1
        \\
        &=\sqrt{1-t}-\frac{1}{2}t\ln\left(\frac{1-\sqrt{1-t}}{\sqrt{1-t}+1}\right)-1+t
        \\
        &\le C_1 t^\nu,
    \end{split}
\end{equation*}
where $C_1$ is a suitable positive constant and $0<\nu<1$\footnote{The computation of the above primitive was performed using Mathematica.}. 
\\
Now, by using the change of variable $z=x-1+t$, we estimate $T_2$ as follows
\begin{equation*}
    \begin{split}
        T_2(t)&= \int_{1-t}^1\varphi\left(\left|f(x-1+t)-f(x)\right|\right) dx
        \\
        &=\int_0^t\varphi\left(\left|f(z)-f(z-t+1)\right|\right) dz 
    	 \\
        &=\int_0^t\left(\sqrt{1+\frac{1-t}{z}}-1\right)dz 
        \\
        &=\lim_{y\to 0}\left[(-1 + \sqrt{-(-1 + t - z)/z}) z + \frac{\sqrt{1 - t} \sqrt{-(-1 + t - z)/z} \sqrt{z} \sinh^{-1}(\frac{\sqrt{z}}{\sqrt{1 - t}})}{\sqrt{\frac{-1 + t - z}{-1 + t}}}\right]_y^t
        \\
        &\le C_2 t^{\frac{1}{2}},
    \end{split}
\end{equation*}
where $C_2$ is a suitable positive constant\footnote{The computations were again performed using Mathematica. In particular, to determinate the exact order one can consider the Taylor expansion of $T_2(t)$ at $t=0$.}. Hence, $f\in \widetilde{Lip}(\nu,\varphi)$ for every $0<\nu\leq\frac{1}{2}$.
\\
Now, we take $\lambda=2$. Arguing as before, by the change of variable $x+t=z$ we get
\begin{equation*}
    \begin{split}
        \int_\Omega\varphi\left(2\left|f(x+t)-f(x)\right|\right) dx&\geq\int_0^{1-t}\varphi\left(2\left|f(x+t)-f(x)\right|\right)dx
        \\
        &=\int_t^1\varphi\left(2\left|f(z)-f(z-t)\right|\right) dz
        \\
        &=\int_t^1\frac{t}{z-t} dz=+\infty.
    \end{split}
\end{equation*}
Therefore, by Proposition \ref{noLip}, we conclude that $f\notin Lip(\nu,\varphi)$, with $0<\nu\le 1$.
\\
In summary, $f\in \widetilde{Lip}(\nu,\varphi)\setminus Lip(\nu,\varphi)$, for every $0<\nu\leq \frac{1}{2}$, and this proves that, in general, the inclusion stated in Proposition \ref{inclusion} is proper.
\end{example}
In the following example, we provide a case when the inclusion given in Proposition \ref{inclusion} is actually an equality.
\begin{example}  
\rm
    We consider $\varphi(u)=u^p$ with $1\le p<+\infty$. In this case, it is well-known that $I^\varphi[\,\cdot\,]=\|\cdot\|_p^p$ and then $Lip(\nu,\varphi) \subset \widetilde{Lip}(\nu,\varphi) \subset Lip(\nu/p,p)$, with $0<\nu\le 1$ (see \cite{EstimatesDurr}). Obviously, if $p=1$, we get:
$\widetilde{Lip}(\nu,\varphi) = Lip(\nu,1)$.
\end{example}
As a consequence of Theorem \ref{mainthm3orlicz}, we can state the following.
 \begin{corollary} Let $f\in \widetilde{Lip}(\nu,\varphi)$, with $0<\nu\le1$ and $\varphi$ be a convex $\varphi$-function. Thus, under the assumptions made in Theorem \ref{thm3}, there exist $\lambda>0$ and $\widetilde{\Gamma}_1>0$ such that
\[
I^\varphi\left[\lambda\left( K_nf-f\right)\right]\le \widetilde{\Gamma}_1 n^{-\nu/2},
\]
for sufficiently large $n\in\N$, with $\widetilde{\Gamma}_1:=\Lambda_1 \widetilde{C}$ where $\Lambda_1$ and $\widetilde{C}$ are suitable constants arising from Theorem \ref{mainthm3orlicz} and (\ref{Lipphiweak}), respectively. 
\end{corollary}

This qualitative result can be enhanced under the assumptions of Theorem \ref{secondmodestimate}. Specifically, in this context, we achieve a better rate of modular convergence, as it is shown by the following.
\begin{corollary} Let $f\in \widetilde{Lip}(\nu,\varphi)$, with $0<\nu\le1$ and $\varphi$ be a convex $\varphi$-function. Thus, under the assumptions made in Theorem \ref{secondmodestimate}, there exist $\lambda>0$ and $\widetilde{\Gamma}_2>0$ such that
\[
I^\varphi\left[\lambda\left( K_nf-f\right)\right]\le \widetilde{\Gamma}_2 n^{-\nu},
\]
for sufficiently large $n\in\N$, with $\widetilde{\Gamma}_2:=\Lambda_2 \widetilde{C}$ where $\Lambda_2$ and $\widetilde{C}$ are suitable constants arising from Theorem \ref{secondmodestimate} and (\ref{Lipphiweak}), respectively. 
\end{corollary}
\section{Inverse approximation and characterization of Lipschitz classes} \label{sec7} 

Here, we establish an inverse approximation theorem for functions in Orlicz spaces by means of a suitable $K$-functional and its equivalence with the strong $\varphi$-modulus of smoothness. Combined with the qualitative results presented in Section~\ref{s6}, this provides a complete characterization of the Lipschitz classes $Lip(\nu,\varphi)$, with $0 < \nu < 1$, introduced in~(\ref{Lipphi}), in terms of the convergence rate of the Kantorovich NN operators. 
\vskip0.2cm

  To derive the inverse approximation theorem, we first establish two norm estimates for the derivatives of the operators. We start with the following Bernstein-type inequality in the Luxemburg norm.
\begin{theorem}\label{bernstein-type-inequality}
    Let $\varphi$ be a convex $\varphi$-function, and $\sigma$ be a sigmoidal function such that the resulting $\phis$ satisfies $\phis^\prime\in L^1(\R)$ and $M_0(\phis^\prime)<+\infty$. Thus, for every $f\in L^\varphi(I)$ there holds
    \be\label{norm-operator-estimate1}
    \left\|K_n^\prime f\right\|_\varphi\leq C_\sigma n \left\|f\right\|_\varphi,
    \ee
where $C_\sigma>0$ is a suitable absolute constant depending only on $\sigma$.
\end{theorem}
\begin{proof} Let $\lambda>0$ be fixed.
    From $(S 2)$, $\phis$ and also $K_nf$ are differentiable. Thus, we can write
    \begin{equation}\label{derivativeop}
        \begin{split}    
    \left(K^\prime_nf\right)(x)&=n\frac{\displaystyle\sum_{k=\lceil na \rceil}^{\lfloor nb \rfloor-1}\phis^\prime(nx-k)n\int_{k/n}^{(k+1)/n}f(u)du}{\displaystyle\sum_{i=\lceil na \rceil}^{\lfloor nb \rfloor-1} \phis(nx-i)}\
    \\
    &\qquad -n\frac{\displaystyle\left(\sum_{k=\lceil na \rceil}^{\lfloor nb \rfloor-1}\phis(nx-k)n\int_{k/n}^{(k+1)/n}f(u)du\right) \left(\sum_{i=\lceil na \rceil}^{\lfloor nb \rfloor-1} \phis^\prime(nx-i)\right) }{\left(\displaystyle\sum_{i=\lceil na \rceil}^{\lfloor nb \rfloor-1}\phis(nx-i)\right)^2},
    \end{split}
    \end{equation}
    $x\in I$. By using the convexity of the $\varphi$ and (\ref{a3}), we can estimate the derivative of the operators as follows
    \begin{equation}\label{firststep}
        \begin{split}   &I^\varphi\left[\lambda\left(K^\prime_nf\right)\right]
        \\
        &\leq \int_a^b\varphi\left(\frac{2\lambda n}{\phis(2)}\sum_{k=\lceil na \rceil}^{\lfloor nb \rfloor-1}\left|\phis^\prime(nx-k)\right|n\int_{k/n}^{(k+1)/n}\left|f(u) \right| \ du \right ) \ dx
        \\
        &\qquad+\int_a^b\varphi\left(\frac{2\lambda n}{\phis^2(2)}\sum_{k=\lceil na \rceil}^{\lfloor nb \rfloor-1}\phis(nx-k)\sum_{i=\lceil na \rceil}^{\lfloor nb \rfloor-1}\left|\phis^\prime(nx-i)\right| n\int_{k/n}^{(k+1)/n}\left|f(u) \right| \ du \right )\ dx
        \\
        &=:I_1+I_2.
        \end{split}
    \end{equation}
    As for $I_1$, applying Jensen inequality twice, we get
    \begin{equation*}
        \begin{split}
            I_1&=\int_a^b\varphi\left(\frac{2\lambda n}{\phis(2)}\sum_{k=\lceil na \rceil}^{\lfloor nb \rfloor-1}\left|\phis^\prime(nx-k)\right|n\int_{k/n}^{(k+1)/n}\left|f(u) \right| \ du \right ) \ dx
            \\
            &\leq \frac{1}{M_0(\phis^\prime)} \int_a^b\sum_{k=\lceil na \rceil}^{\lfloor nb \rfloor-1}\left|\phis^\prime(nx-k)\right|\varphi\left(\frac{2 M_0(\phis^\prime)\lambda n^2}{\phis(2)}\int_{k/n}^{(k+1)/n}\left|f(u) \right| \ du \right ) \ dx
            \\
            &\leq \frac{1}{M_0(\phis^\prime)} \sum_{k=\lceil na \rceil}^{\lfloor nb \rfloor-1}\int_a^b\left|\phis^\prime(nx-k)\right|\ dx\left[n\int_{k/n}^{(k+1)/n}\varphi\left(\frac{2 M_0(\phis^\prime)\lambda n}{\phis(2)}\left|f(u) \right| \right )  du \right]
            \\
            &\leq \frac{\left\|\phis^\prime\right\|_1}{M_0(\phis^\prime)}\sum_{k=\lceil na \rceil}^{\lfloor nb \rfloor-1}\left[\int_{k/n}^{(k+1)/n}\varphi\left(\frac{2 M_0(\phis^\prime)\lambda n}{\phis(2)}\left|f(u) \right| \right )  du \right]
            \\
            &\leq\frac{\left\|\phis^\prime\right\|_1}{M_0(\phis^\prime)}\int_a^b \varphi\left(\frac{2 M_0(\phis^\prime)\lambda n}{\phis(2)}\left|f(u) \right| \right )  du
            \\
            &= \frac{\left\|\phis^\prime\right\|_1}{M_0(\phis^\prime)} I^\varphi\left[\frac{2 M_0(\phis^\prime)\lambda n}{\phis(2)} f\right],
        \end{split}
    \end{equation*}
   where $$\int_a^b\left|\phis^\prime(nx-k)\right|dx\leq\frac{1}{n}\int_\R \left|\phis^\prime(t)\right|dt=\frac{\left\|\phis^\prime\right\|_1}{n}.$$
   As for $I_2$, we again apply Jensen inequality twice, together with (\ref{m0}), to obtain
   \begin{equation*}
       \begin{split}
           I_2&=\int_a^b\varphi\left(\frac{2\lambda n}{\phis^2(2)}\sum_{k=\lceil na \rceil}^{\lfloor nb \rfloor-1}\phis(nx-k)\sum_{i=\lceil na \rceil}^{\lfloor nb \rfloor-1}\phis^\prime(nx-i)n\int_{k/n}^{(k+1)/n}\left|f(u) \right| \ du \right ) \ dx
           \\
           &\leq \int_a^b\varphi\left(\frac{2 M_0(\phis^\prime)\lambda n}{\phis^2(2)}\sum_{k=\lceil na \rceil}^{\lfloor nb \rfloor-1}\phis(nx-k)n\int_{k/n}^{(k+1)/n}\left|f(u) \right| \ du \right ) \ dx
           \\
           &\leq \int_a^b\sum_{k=\lceil na \rceil}^{\lfloor nb \rfloor-1}\phis(nx-k)\varphi\left(\frac{2 M_0(\phis^\prime)\lambda n^2}{\phis^2(2)}\int_{k/n}^{(k+1)/n}\left|f(u) \right| \ du \right ) \ dx
           \\
           &\leq \sum_{k=\lceil na \rceil}^{\lfloor nb \rfloor-1}\int_a^b\phis(nx-k)\ dx\left[n\int_{k/n}^{(k+1)/n}\varphi\left(\frac{2 M_0(\phis^\prime)\lambda n}{\phis^2(2)}\left|f(u) \right| \right )\ du \right]  
           \\
           &\leq \sum_{k=\lceil na \rceil}^{\lfloor nb \rfloor-1}\int_{k/n}^{(k+1)/n}\varphi\left(\frac{2 M_0(\phis^\prime)\lambda n}{\phis^2(2)}\left|f(u) \right| \right )\ du
           \\
           &\leq  I^\varphi\left[\frac{2 M_0(\phis^\prime)\lambda n}{\phis^2(2)} f\right],
       \end{split}
   \end{equation*}
    where
    $$ \int_a^b\phis(nx-k)\ dx\leq \frac{1}{n}\int_\R \phis(t) dt=\frac{\left\|\phis\right\|_1}{n}=\frac{1}{n},$$
      since $\phi_{\sigma} \in L^1(\R)$ (see Lemma~\ref{lemma1}).\\
    Setting now $C:=2\max\left\{\frac{\left\|\phis^\prime\right\|_1}{M_0(\phis^\prime)},1\right\}$, by the monotonicity of the modular we can write
    \begin{equation}\label{modular-bernstein-type}
    I^\varphi\left[\lambda \left(K_n^\prime f\right)\right]\leq C I^\varphi\left[\frac{2 M_0(\phis^\prime)\lambda n}{\phis^2(2)} f\right].
    \end{equation}
   Since inequality~\eqref{modular-bernstein-type} holds for any fixed \(\lambda > 0\), we can also derive its norm version by proceeding as in the proof of Theorem~\ref{newthm}. Hence, the thesis follows by a suitable positive constant $C_\sigma$.
\end{proof}
\begin{theorem}
    Let $\varphi$ be a convex $\varphi$-function such that $u\mapsto u^{-\beta}\varphi(u)$, $u>0$, is increasing, with some $\beta>1$. Moreover, let $\sigma$ be a sigmoidal function such that the resulting $\phis$ satisfies $M_1(\phis)+M_1(\phis^\prime)<+\infty$. Thus, for every $f\in W^{1,\varphi}(I)$ there holds
    \be\label{norm-operator-estimate2}
    \left\|K_n^\prime f\right\|_\varphi\leq D_\sigma  \left\|f^\prime\right\|_\varphi,
    \ee
where $D_\sigma>0$ is a suitable absolute constant depending only on $\sigma$.
\end{theorem}
\begin{proof}
    Let $\lambda>0$. If $\mathbf{1}$ denotes the unitary constant function, then $K_n \mathbf{1} \equiv \mathbf{1}$, 
and hence $K_n^\prime \mathbf{1} \equiv 0$. 
Using this, together with the expression of $K^\prime_nf$ in (\ref{derivativeop}), (\ref{taylor1}) and proceeding as in (\ref{firststep}), we get
    \begin{equation*}
    \begin{split}
        I^\varphi\left[\lambda\left(K_n^\prime f\right)\right]&=I^\varphi\left[\lambda\left(K_n^\prime f-fK_n^\prime \bf{1}\right)\right]
        \\
        &\leq \int_a^b\varphi\left(\frac{\lambda n^2}{\phis(2)}\sum_{k=\lceil na \rceil}^{\lfloor nb \rfloor-1} \left|\phis^\prime(nx-k)\right|\int_{k/n}^{(k+1)/n}\left|\int_x^u\left|f^\prime(t)\right|dt\right|\ du\right.
        \\
       &\qquad \left.+\frac{\lambda n^2}{\phis^2(2)}\sum_{k=\lceil na \rceil}^{\lfloor nb \rfloor-1}\phis(nx-k)\sum_{i=\lceil na \rceil}^{\lfloor nb \rfloor-1}\left|\phis^\prime(nx-i)\right|\int_{k/n}^{(k+1)/n}\left|\int_x^u\left|f^\prime(t)\right|dt\right|\ du\right)\ dx
       \\
       &\leq \int_a^b\varphi\left(\frac{2\lambda n^2}{\phis(2)}\sum_{k=\lceil na \rceil}^{\lfloor nb \rfloor-1}\left|\phis^\prime(nx-k)\right|\int_{k/n}^{(k+1)/n}|u-x|\mathcal{M}f^\prime(x)\ du\right)\ dx
       \\
       &\qquad + \int_a^b\varphi\left(\frac{2\lambda M_0(\phis^\prime) n^2}{\phis^2(2)}\sum_{k=\lceil na \rceil}^{\lfloor nb \rfloor-1}\phis(nx-k)\int_{k/n}^{(k+1)/n}|u-x|\mathcal{M}f^\prime(x)\ du\right)\ dx
       \\
       &=:J_1+J_2.
    \end{split}
    \end{equation*}
As for $J_1$, we proceed with the following estimate
\begin{equation*}
    \begin{split}
        J_1&=\int_a^b\varphi\left(\frac{2\lambda n^2}{\phis(2)}\mathcal{M}f^\prime(x)\sum_{k=\lceil na \rceil}^{\lfloor nb \rfloor-1}\left|\phis^\prime(nx-k)\right|\int_{k/n}^{(k+1)/n}\left|u-\frac{k}{n}+\frac{k}{n}-x\right|\ du\right)\ dx
        \\
        &\leq \int_a^b\varphi\left(\frac{2\lambda n}{\phis(2)}\mathcal{M}f^\prime(x)\sum_{k=\lceil na \rceil}^{\lfloor nb \rfloor-1}\left|\phis^\prime(nx-k)\right|\int_{k/n}^{(k+1)/n}\left(\left|nu-k\right|+\left|k-nx\right|\right)\ du\right)\ dx
        \\
        &\leq \int_a^b\varphi\left(\frac{2\lambda n}{\phis(2)}\mathcal{M}f^\prime(x)\sum_{k=\lceil na \rceil}^{\lfloor nb \rfloor-1}\left|\phis^\prime(nx-k)\right|\int_{k/n}^{(k+1)/n}\left(1+\left|k-nx\right|\right)\ du\right)\ dx
        \\
        &= \int_a^b\varphi\left(\frac{2\lambda }{\phis(2)}\mathcal{M}f^\prime(x)\sum_{k=\lceil na \rceil}^{\lfloor nb \rfloor-1}\left|\phis^\prime(nx-k)\right|\left(1+\left|k-nx\right|\right)\right)\ dx
        \\
        &\leq \int_a^b\varphi\left(\frac{2\lambda }{\phis(2)}\left(M_0(\phis^\prime)+M_1(\phis^\prime)\right)\mathcal{M}f^\prime(x)\right)\ dx
        \\
        &=I^\varphi\left[\frac{2\lambda}{\phis(2)}\left(M_0(\phis^\prime)+M_1(\phis^\prime)\right)\mathcal{M}f^\prime\right],
    \end{split}
\end{equation*}
where $M_0(\phis^\prime)<+\infty$ since $M_1(\phis^\prime)<+\infty$. Similarly, for $J_2$ we have
\begin{equation*}
    \begin{split}
        J_2&\leq \int_a^b\varphi\left(\frac{2M_0(\phis^\prime)\lambda n}{\phis^2(2)}\sum_{k=\lceil na \rceil}^{\lfloor nb \rfloor-1}\phis(nx-k)\int_{k/n}^{(k+1)/n}\left(\left|nu-k\right|+\left|k-nx\right|\right)\mathcal{M}f^\prime(x)\ du\right)\ dx
        \\
        &\leq \int_a^b\varphi\left(\frac{2M_0(\phis^\prime)\lambda }{\phis^2(2)}\mathcal{M}f^\prime(x)\sum_{k=\lceil na \rceil}^{\lfloor nb \rfloor-1}\phis(nx-k)\left(1+\left|k-nx\right|\right)\right)\ dx
        \end{split}
\end{equation*}
\begin{equation*}
    \begin{split}
        &\leq \int_a^b\varphi\left(\frac{2M_0(\phis^\prime)\lambda }{\phis^2(2)}\left(M_0(\phis)+M_1(\phis)\right)\mathcal{M}f^\prime(x)\right)\ dx
        \\
        &=I^\varphi\left[\frac{2M_0(\phis^\prime)\lambda }{\phis^2(2)}\left(1+M_1(\phis)\right)\mathcal{M}f^\prime\right].
    \end{split}
\end{equation*}
Proceeding as in the proof of Theorem~\ref{bernstein-type-inequality}, the above inequalities can be combined into a single one, which is valid for any fixed \(\lambda > 0\), while arguing as in the proof of Theorem~\ref{newthm}, we then derive the corresponding Luxemburg norm estimate.
Finally, by exploiting Theorem~\ref{Mbounded}, we can get the thesis.
\end{proof}
Now, to obtain an inverse approximation result in $L^\varphi(I)$, 
we need to use a specific tool, namely the so-called \textit{$K$-functional}, recalled below.
\\
Given $f\in L^\varphi(I)$, we have
\be\label{K-funct}
\mathcal{K}(f,\delta)_\varphi:=\inf_{g\in W^{1,\varphi}(I)}\left\{\left\|f-g\right\|_\varphi+\delta \left\|g^\prime\right\|_\varphi \right\}.
\ee
As done in \cite{Garidi1991}, it can be proved that $K(f,\delta)_\varphi$ and $\omega(f,\delta)_\varphi$ are equivalent, in the sense that there exist $c_1,c_2>0$ such that
\be\label{eqKfunct}
c_1 \omega(f,\delta)_\varphi\leq \mathcal{K}(f,\delta)_\varphi \leq c_2 \omega(f,\delta)_\varphi.
\ee
Moreover, we also recall the lemma established by Berens and Lorentz in \cite{BerensLorentz1972}.
\begin{lemma}[\cite{BerensLorentz1972}]\label{berensLorentz}
    For any fixed $r$ and $\nu$ such that $r>\nu$, we suppose that the positive sequence $\eta(n)$, $n\in\N$, such that $$\eta(n)\leq A k^{-\nu}+B(k/n)^r\eta(k), \ \eta(n_0)\leq C,$$ for all $n\geq n_0$ and $k$, where $A,B, C>0$. Thus, there exists $D>0$ such that
    $$\eta(n)\leq D n^{-\nu},$$
    for all $n\geq n_0$.
\end{lemma}
Now, we are ready to provide the following inverse approximation theorem.
\begin{theorem}\label{inversionthm}
    Let $\varphi$ be a convex $\varphi$-function such that $u\mapsto u^{-\beta}\varphi(u)$, $u>0$, is increasing, with some $\beta>1$. Moreover, let $\sigma$ be a sigmoidal function such that the resulting $\phis$ satisfies $\phis^\prime\in L^1(\R)$ and $M_1(\phis)+M_1(\phis^\prime)<+\infty$. Thus, for every $f\in L^\varphi(I)$ and $0<\nu<1$ such that
    \begin{equation}\label{inversion-eq}
    \left\|K_nf-f\right\|_\varphi=\mathcal{O}(n^{-\nu}),\qquad n\to+\infty,
    \end{equation}
    it turns out that
    \[
    \omega(f,\delta)_\varphi=\mathcal{O}(\delta^{\nu}),\qquad \delta\to 0^+.
    \]
\end{theorem}
\begin{proof}
    Let $0<\nu<1$ and $n\in\N$. From (\ref{inversion-eq}), there exist a constant $A>0$ and $\overline{k}\in\N$ such that $\left\|K_kf-f\right\|_\varphi\leq A k^{-\nu}$, for every $k\geq \overline{k}$. Let now $k\geq \overline{k}$. Taking into account (\ref{K-funct}) with $\delta=1/n$, along with the norm estimates (\ref{norm-operator-estimate1}) and (\ref{norm-operator-estimate2}), there exists a suitable $B>0$ such that
    \begin{equation*}
        \begin{split}
            \mathcal{K}(f,1/n)_\varphi&\leq \left\|f-K_kf\right\|_\varphi+\frac{1}{n}\left\|K_k^\prime f\right\|_\varphi
            \\
            &\leq  \left\|f-K_kf\right\|_\varphi+\frac{1}{n}\left(\left\|K_k^\prime (f-g)\right\|_\varphi+\left\|K_k^\prime g\right\|_\varphi\right)
            \\
            &\leq  \left\|f-K_kf\right\|_\varphi+\frac{Bk}{n}\left(\left\|f-g\right\|_\varphi+\frac{1}{k}\left\| g^\prime\right\|_\varphi\right),
        \end{split}
    \end{equation*}
    where $g\in W^{1,\varphi}(I)$. Now, passing to the infimum with respect to $g$ and using (\ref{inversion-eq}), we have
    \begin{equation*}
        \begin{split}
            \mathcal{K}(f,1/n)_\varphi&\leq \left\|f-K_kf\right\|_\varphi+\frac{Bk}{n}\mathcal{K}(f,1/k)_\varphi
            \\
            &\leq Ak^{-\nu}+\frac{Bk}{n}\mathcal{K}(f,1/k)_\varphi.
        \end{split}
    \end{equation*}
    Taking $\eta(n):=\mathcal{K}(f,1/n)_\varphi$ and $r=1$ in Lemma \ref{berensLorentz}, and using the equivalence in (\ref{eqKfunct}), we obtain that
    \[
    \omega(f,1/n)_\varphi=\mathcal{O}(n^{-\nu}),\qquad0<\nu<1,
    \]
    for sufficiently large $n\geq n_0\in\N$. In particular, for any given $0<\delta<\frac{1}{n_0}$, there exists a suitable $n\geq n_0$ such that $\frac{1}{2n}\leq \delta\leq\frac{1}{n}$. Therefore, by the monotonicity of the modulus of smoothness, we conclude that
    \[
    \omega(f,\delta)_\varphi\leq \omega(f,1/n)_\varphi=\mathcal{O}(n^{-\nu})=\mathcal{O}(\delta^\nu),\qquad\delta\to0^+.
    \]
    This completes the proof.
\end{proof}
\begin{example}\rm
    As an example of a sigmoidal activation function satisfying the assumptions required in Theorem \ref{inversionthm}, 
we can mention the logistic function $\sigma_l$ defined in (\ref{logisticf}), whose density function is given by
\[
\phi_{\sigma_l}(x) = \frac{1}{2}\big(\sigma_l(x+1) - \sigma_l(x-1)\big)
= \frac{e^2 - 1}{2} \cdot \frac{1}{(1 + e^{1+x})(1 + e^{1-x})}, \qquad x \in \mathbb{R}.
\]
It is easy to see that
\[
\phi_{\sigma_l}'(x) = -\frac{e(e^2 - 1)}{2} \cdot 
\frac{e^x - e^{-x}}{(1 + e^{1+x})^2 (1 + e^{1-x})^2}, \qquad x \in \mathbb{R}.
\]
Hence, due to its exponential decay, we can conclude that $\phi_{\sigma_l} \in C^1(I) \subset W^{1,1}(I)$, 
and the moment-type conditions 
$M_1(\phi_{\sigma_l}) +  M_1(\phi_{\sigma_l}^\prime) < +\infty$ 
are also satisfied. Clearly, all these conditions are trivially fulfilled for compactly supported density functions.
\end{example}
Now, by combining the direct qualitative result given in Corollary \ref{cordirectresult} 
with the inverse approximation theorem presented above (Theorem \ref{inversionthm}), 
we obtain the following final result, which provides a complete characterization 
of the Orlicz-Lipschitz classes $Lip(\nu,\varphi)$ for every $0 < \nu < 1$.
\begin{theorem}\label{characterizationLip}
  Let $\varphi$ be a convex $\varphi$-function such that $u\mapsto u^{-\beta}\varphi(u)$, $u>0$, is increasing, with some $\beta>1$ and $\sigma$ be such that the resulting $\phis$ satisfies $\phis^\prime\in L^1(\R)$. If $\varphi$ is also an $N$-function (risp. satisfying the $\Delta^\prime$-condition) and $\phis$ is compactly supported (risp. satisfying $M_{0,1}^\varphi(\phis)+M_1(\phis)+M_1(\phis^\prime)<+\infty$), then $$f\in Lip(\nu,\varphi),\qquad0<\nu<1,$$ if and only if
  \[
    \left\|K_nf-f\right\|_\varphi=\mathcal{O}(n^{-\nu}),\qquad n\to+\infty.
    \]
\end{theorem}

\section{Final remarks, conclusions and open problems} \label{sec8} 

In this work, we provide a comprehensive asymptotic analysis for a family of NN operators based on both strong (i.e., with respect to the Luxemburg norm) and weak (modular) sharp bounds in Orlicz spaces, using suitable versions of the $\varphi$-modulus of smoothness.\\
These results are derived using density theorems and asymptotic estimates in the context of Sobolev-Orlicz spaces. To this aim, we employ new tools and techniques of proof, that are summarized below.
\begin{enumerate}
    \item The Orlicz Minkowski inequality, which allows to obtain unifying Luxemburg norm-based estimates for the operators, suitable for several instances of Sobolev-Orlicz and Orlicz spaces.
    \item A new weak version of the aforementioned Orlicz Minkowski inequality, which holds under less restrictive assumptions on the $\varphi$-function with respect to its strong (classical) counterpart.
    \item The introduction of a new notion of discrete absolute $\varphi$-moments of the hybrid type (see (\ref{phi-moment})).
    \item The application of the HL-maximal function in Orlicz spaces.
    \item The introduction of the new space $\mathcal{W}^{1,\varphi}(I)$ (see (\ref{subspace})), embedded in $ W^{1,\varphi}(I)$ and modularly dense in $L^\varphi(I)$, which allows to derive an asymptotic estimate under more general assumptions (with respect to the previous results) on the involved $\varphi$-function.
    \item A new weak (modular) extension of the Sobolev-Orlicz density result given by H. Musielak via Steklov functions (Lemma \ref{lemmaweak}).
    \item The study of the inclusion properties between weak and strong Orlicz Lipschitz classes, defined by different moduli of smoothness.
    \item A characterization theorem for the Lipschitz class $Lip(\nu,\varphi)$, with $0<\nu<1$, by the rate of convergence of Kantorovich NN operators.
\end{enumerate}
In particular, the result mentioned at point 8. proves the sharpness of the estimates given in Corollary \ref{cordirectresult} for $0<\nu<1$.
As an open problem, it remains to consider the case corresponding to the saturation order, i.e., for $\nu=1$. In general, saturation theorems are highly delicate and require new strategies of proof.
Moreover, the weak Lipschitz classes introduced here, $\widetilde{Lip}(\nu,\varphi)$, could also be characterized in terms of the rate of convergence. To do this, and in order to adopt a strategy similar to the one used in Theorem \ref{inversionthm}, a suitable version of the $K$-functional equivalent to the weak $\varphi$-modulus of smoothness should be introduced and studied.
\vskip0.2cm

 Concluding these final remarks, we want to stress that the study of approximation results in Orlicz and Sobolev-Orlicz spaces can be interesting from the point of view of Approximation Theory since, in these contexts, we are in more general situations than the classical $L^p$-spaces of the classical Sobolev spaces $W^{1, p}$.
\vskip0.2cm

  Indeed, from the general theory of Orlicz and Sobolev-Orlicz spaces that, if the $\varphi$-function is such that $\varphi(t)/t^p$, $1 \miu p <+\infty$, is {\em almost increasing}, i.e., there exists $a \mau 1$ for which:
\be \label{aInc-p}
{\varphi(s) \over s^p } \miu\ a\, {\varphi(t) \over t^p }, 
\ee
for all $0<s<t$, then, for a given compact $I \subset \R$ there holds the following embedding (see Lemma 6.1.6 of \cite{harjulehto2019orlicz}):
\be \label{emb}
W^{1, \varphi}(I)\ \hookrightarrow\ W^{1,p}(I).
\ee
It is well-known that any convex $\varphi$-function $\varphi$ satisfies (\ref{aInc-p}) with $p=1$ (see eqn. (2.2) of \cite{alberico2024modulus}) but not in general for $p>1$, and consequently it is possible to see that the embedding (\ref{emb}), for $p>1$, is not always true. 
Therefore, we can deduce that Sobolev-Orlicz spaces are more general than the classical Sobolev spaces, as well as Orlicz spaces are more general than $L^p$-spaces.
\vskip0.2cm
Indeed, consider for instance the following convex $\varphi$-function (see \cite{rao1991theory}, p. 29):
\be
\widetilde{\varphi}(t)\ :=\ {t^p \over \log(e + t)}, \quad t \mau 0, \quad 1 < p <+\infty,
\ee
and the function $f:[0,1] \to \R$ defined by:
\be
f(x)\ :=\ \begin{cases}
\left({1 \over x\, \ln(1/x)}\right)^{1/p}, \quad x \in (0,1/2) \\
\\
0, \quad \quad \quad \quad \quad \quad \quad \mbox{otherwise}.
\end{cases}
\ee
Now, posed:
\be
u(x)\ :=\ \int_0^x f(t)\, dt, \quad x \in [0,1],
\ee 
it turns out that $u \in W^{1, \widetilde{\varphi}}([0,1])$ while $u \notin W^{1,p}([0,1])$. Indeed, $I^{\widetilde{\varphi}}[u] < +\infty$ and $\|u\|_p<+\infty$ since $u \in AC([0,1])$, and:
$$
I^{\widetilde{\varphi}}[u']\ =\ \int_0^{1/2} \widetilde{\varphi}\left(  f(x) \right)\, dx\ \miu\ C\,  \int_0^{1/2} {1 \over x \ln^2(1/x)}\, dx\ =\ C\, \int_{\ln 2}^{+\infty} {1 \over t^2}\, dt\ <\ +\infty,
$$
since:
$$
\widetilde{\varphi}(f(x))\ =\ {\cal O}\left({1 \over x \ln^2(1/x)} \right), \quad as\ \quad x \to 0^+,
$$
while:
$$
\|u'\|^p_p\ =\ \|f\|^p_p\ =\ \int_0^{1/2} {1 \over x\, \ln(1/x)}\, dx\ =\ \int_{\ln 2}^{+\infty} {1 \over t}\, dt\ =\ +\infty.
$$
As a consequence, we also have that $f \in L^{\widetilde{\varphi}}([0,1])$ and $f \notin L^p([0,1])$. These examples show the consistency of the study of general approximation results in the framework of Orlicz and Sobolev-Orlicz spaces.  
\vskip0.2cm

 Furthermore, we want to stress that Orlicz spaces, and the asymptotic analysis developed in Sobolev-Orlicz spaces, due to the flexibility in the choice of $\varphi$, are useful in certain practical situations where classical $L^p$-approximation fails.
\vskip0.1cm
  For instance, Orlicz approximation allows a suitable analysis for functions having localized peaks which model impulsive noise (see, e.g., \cite{Nikolova2004} for a variational-type approach). In this regard, if we consider a function $f_\varepsilon(x) = 1/\varepsilon$ on a very small interval $[0,\varepsilon]$ and 0 elsewhere, the computation of its $L^p$-norm shows that such 
approximation process is strongly unstable, in the sense that
\begin{equation*}
\|f_\varepsilon\|_p^p=\int_0^\varepsilon \varepsilon^{-p} \ dx=\varepsilon^{1-p}\to +\infty,\qquad\text{as $\varepsilon\to 0^+$.}
\end{equation*}
While, using, e.g., the Orlicz modular $I^\varphi$ generated by a Zygmund-type function $\varphi(t) = t\log(1+t)$, we can reduce the effect of the peak, giving a more balanced approximation (reducing the instability), namely
\begin{equation*}
I^\varphi\left[f_\varepsilon\right]=\log\left(1+\frac{1}{\varepsilon}\right),
\end{equation*}
which clearly increases more slowly than $\varepsilon^{1-p}$ as $\varepsilon\to0^+$. A similar situation occurs in polynomial approximation problems involving functions with local peaks, when the global approximation error must be estimated (see, e.g., \cite{Plesniak1985,Garidi1991}).
\vskip0.2cm
  Moreover, Orlicz approximation is particularly useful for functions having super-polynomial or exponential growth, which do not belong to any $L^p$ space (see \cite{BoscagginColasuonnoNorisSani}). 
This provides a suitable functional framework to analyze and approximate solutions of differential equations with exponential or super-polynomial nonlinearities (see, e.g., \cite{BrezisMerle1991,Ioku}). 
\vskip0.2cm

Finally, an interesting direction for future developments would be the extension of the theoretical results here presented to the multidimensional setting, using the multivariate version of $K_nf$ introduced in \cite{CS14}. Here, several technical difficulties may arise; for instance, the density results via Steklov functions should be extended for functions of several variables, and so on. In any case, the present paper, and its multivariate extension, could open the way to the study of further applications.

\section*{Acknowledgments}

{\small 
The authors would like to thank the anonymous referees for their comments which were very useful in improving the quality of the present manuscript.
\\
The authors are members of the Gruppo Nazionale per l'Analisi Matematica, la Probabilit\`a e le loro Applicazioni (GNAMPA) of the Istituto Nazionale di Alta Matematica (INdAM), of the network RITA (Research ITalian network on Approximation), and of the UMI (Unione Matematica Italiana) group T.A.A. (Teoria dell'Approssimazione e Applicazioni). 
}

\section*{Funding}

{\small Both the authors have been supported by: PRIN 2022 PNRR: ``RETINA: REmote sensing daTa INversion with multivariate functional modeling for essential climAte variables characterization", funded by the European Union under the Italian National Recovery and Resilience Plan (NRRP) of NextGenerationEU, under the Italian Minister of University and Research MUR (Project Code: P20229SH29, CUP: J53D23015950001).
}

\section*{Conflicts of interests/Competing interests}
{\small The authors declare that they have not conflict of interest and/or competing interest.}

\section*{Availability of data and material and Code availability}

{ \small Not applicable.}



\begin{thebibliography}{99}


\bibitem{alberico2024modulus} A. Alberico, A. Cianchi, L. Pick, L. Slavíkova, On the modulus of continuity of fractional Orlicz-Sobolev functions, \textit{Math. Ann.}, 391 (2025), 2429--2477.

\bibitem{AN1997} G.A. Anastassiou, Rate of convergence of some neural network operators to the unit-univariate case, {\em J. Math. Anal. Appl.},  212(1) (1997), 237--262.

\bibitem{BMV2003} C. Bardaro, J. Musielak, G. Vinti, Nonlinear integral operators and applications, Walter de Gruyter, 2003.

\bibitem{BAR} A.R. Barron, Universal approximation bounds for superpositions of a sigmoidal function, {\em IEEE Trans. Inform. Theory}, 39(3) (1992), 930-945.

\bibitem{BR80} C. Bennett, K. Rudnick, On Lorentz-Zygmund spaces, \textit{Diss. Math.}, 175 (1980), 1--72.

\bibitem{BerensLorentz1972} H. Berens, G. Lorentz, Inverse theorems for Bernstein polynomials, \textit{Indiana Univ. Math. J.}, 21 (1972), 693--708.

\bibitem{BoscagginColasuonnoNorisSani} A. Boscaggin, F. Colasuonno, B. Noris, F. Sani, An Orlicz space approach to exponential elliptic problems in higher dimensions, preprint, \url{arXiv:2503.16105}, 2025.

\bibitem{BrezisMerle1991}
H. Brezis, F.Merle, Uniform estimates and blow-up behavior for solutions of $-\Delta u = V(x)e^{u}$ in two dimensions, {\em Communications in Partial Differential Equations}, 16 (1991), 1223--1253.

\bibitem{CACO2025} M. Cantarini, D. Costarelli, Simultaneous approximation by neural network operators with applications to Voronovskaja formulas, {\em Math. Nachr.}, 298(3) (2025), 871--885.

\bibitem{CAEU1} P. Cardaliaguet, G. Euvrard, Approximation of a function and its derivative with a neural network, \textit{Neural Netw.}, 5 (1992), 2, 207--220.

\bibitem{CACH1} Z. Chen, F. Cao, The approximation operators with sigmoidal functions, \textit{Computers Math. Appl.}, 58 (2009), 4, 758--765.

\bibitem{CACH2} Z. Chen, F. Cao, The construction and approximation of a class of neural networks operators with ramp functions, \textit{J. Comput. Anal. Appl.}, 14 (2012), 1, 101--112.

\bibitem{chen2022construction} H. Chen, D. Yu, Z. Li, The construction and approximation of ReLU neural network operators, \textit{J. Function Spaces}, 1 (2022), 1713912.

\bibitem{IE1972} I. S. Cheng, J. J. Kozak, Application of the theory of Orlicz spaces to statistical mechanics. I. Integral equations, \textit{J. Math. Phys.}, 13 (1972), 1, 51--58.

\bibitem{CGW21} I. Chlebicka, P. Gwiazda, A. Świerczewska-Gwiazda, A. Wróblewska-Kamińska, Partial differential equations in anisotropic Musielak-Orlicz spaces, Springer, 2021.

\bibitem{CC24} L. Coroianu, D. Costarelli, Best Approximation and Inverse Results for Neural Network Operators, \textit{Results Math.}, 79 (2024), 5, 193.

\bibitem{coroianu2024approximation} L. Coroianu, D. Costarelli, M. Natale, A. Pantiș, The approximation capabilities of Durrmeyer-type neural network operators, \textit{J. Appl. Math. Comput.}, 1--19 (2024).

\bibitem{CO2} D. Costarelli, Density results by deep neural network operators with integer weights, \textit{Math. Model. Anal.}, 27 (2022), 4, 547--560.

\bibitem{CJat23} D. Costarelli, Approximation error for neural network operators by an averaged modulus of smoothness, \textit{J. Approx. Theory}, 294 (2023), 105944.

\bibitem{costarelli2024convergence} D. Costarelli, Convergence and high order of approximation by Steklov sampling operators, \textit{Banach J. Math. Anal.}, 18 (2024), 4, 1--25.

\bibitem{cpNNlp} D. Costarelli, M. Piconi, Asymptotic Analysis of Neural Network Operators Employing the Hardy-Littlewood Maximal Inequality, \textit{Mediterr. J. Math.}, 21 (2024), 7, 1--24.

\bibitem{EstimatesDurr} D. Costarelli, M. Piconi, G. Vinti, Quantitative estimates for Durrmeyer-sampling series in Orlicz spaces, \textit{Sampling Theory, Signal Process. Data Anal.}, 21 (2023), 3.

\bibitem{costarelli2018approximation} D. Costarelli, A. R. Sambucini, Approximation results in Orlicz spaces for sequences of Kantorovich max-product neural network operators, \textit{Results Math.}, 73 (2018), 1, 15.

\bibitem{COSP1} D. Costarelli, R. Spigler, Approximation results for neural network operators activated by sigmoidal functions, \textit{Neural Netw.}, 44 (2013), 101--106.

\bibitem{CS14} D. Costarelli, R. Spigler, Convergence of a family of neural network operators of the Kantorovich type, \textit{J. Approx. Theory}, 185 (2014), 80--90.

\bibitem{COSP2} D. Costarelli, R. Spigler, How sharp is the Jensen inequality?, \textit{J. Inequal. Appl.}, 2015 (2015), 1--10.

\bibitem{CVMathNachr} D. Costarelli, G. Vinti, Convergence for a family of neural network operators in Orlicz spaces, \textit{Math. Nachr.}, 290 (2017), 2-3, 226--235.

\bibitem{CY} G. Cybenko, Approximation by superpositions of a sigmoidal function, \textit{Math. Control Signals Systems}, 2 (1989), 4, 303--314.

\bibitem{DDFH21} I. Daubechies, R. DeVore, S. Foucart, B. Hanin, G. Petrova, Nonlinear approximation and (deep) ReLU networks, \textit{Constructive Approx.}, 55 (2022), 1, 127--172.

\bibitem{DELO1} R. A. DeVore, G. G. Lorentz, Constructive approximation, Springer Science \& Business Media, 1993, vol. 303.

\bibitem{D71} T. Donaldson, Nonlinear elliptic boundary value problems in Orlicz-Sobolev spaces, \textit{J. Differ. Equ.}, 10 (1971), 3, 507--528.

\bibitem{DTOS} T. K. Donaldson, N. S. Trudinger, Orlicz-Sobolev spaces and imbedding theorems, \textit{J. Funct. Anal.}, 8 (1971), 1, 52--75.

\bibitem{DF16} G. Dong, X. Fang, Differential equations of divergence form in separable Musielak-Orlicz-Sobolev spaces, \textit{Boundary Value Probl.}, 2016 (2016), 1--19.

\bibitem{GAO} B. Gao, Y. Xu, Univariant approximation by superpositions of a sigmoidal function, {\em J. Math. Anal. Appl.}, 178(1) (1993) 221--226.

\bibitem{GMS99} M. García-Huidobro, V. Khoi Le, R. Manásevich, K. Schmitt, On principal eigenvalues for quasilinear elliptic differential operators: an Orlicz-Sobolev space setting, \textit{Nonlinear Differ. Equ. Appl.}, 6 (1999), 2, 207--226.

\bibitem{Garidi1991} W. Garidi, On approximation by polynomials in Orlicz spaces, \textit{Approx. Theory Appl.}, 7 (1991), 97--110.

\bibitem{Grafakos2014} L. Grafakos, Maximal Functions, Fourier Transform, and Distributions, in \textit{Classical Fourier Analysis}, Springer New York, 2014, 85--172, ISBN 978-1-4939-1194-3.

\bibitem{MaximalOp} L. Grafakos, D. He, P. Honzík, Maximal operators associated with bilinear multipliers of limited decay, \textit{J. Analyse Math.}, 143 (2021), 231--251.

\bibitem{GR21} I. Gühring, M. Raslan, Approximation rates for neural networks with encodable weights in smoothness spaces, \textit{Neural Netw.}, 134 (2021), 107--130.

\bibitem{harjulehto2019orlicz} P. Harjulehto, P. H\"{a}st\"{o}, Orlicz Spaces and Generalized Orlicz Spaces, \textit{Lect. Notes Math.}, 2236 (2019), 1--12.

\bibitem{Hasto15} P. Hästö, The maximal operator on generalized Orlicz spaces, \textit{J. Funct. Anal.}, 269 (2015), 12, 4038--4048.

\bibitem{H03} S. Hencl, A sharp form of an embedding into exponential and double exponential spaces, \textit{J. Funct. Anal.}, 204 (2003), 1, 196--227.

\bibitem{h1} H. Hudzik, On generalized Orlicz-Sobolev space, \textit{Funct. Approxim. Comment. Math.}, 4 (1976), 37--51.

\bibitem{h2} H. Hudzik, A generalization of Sobolev spaces. I, \textit{Funct. Approxim. Comment. Math.}, 2 (1976), 67--73.

\bibitem{h4} H. Hudzik, The problem of separability, duality, reflexivity and of comparison for generalized Orlicz-Sobolev spaces ${W}_{M}^{k}({\Omega})$, \textit{Comment. Math.}, 21 (1979), 2.

\bibitem{Ioku} N. Ioku, The Cauchy problem for heat equations with exponential nonlinearity, {\em J. Diff. Eq.}, 251 (2011), 1172--1194.

\bibitem{kadak2022neural} U. Kadak, Neural network operators of fuzzy n-cell number valued functions and multidimensional fuzzy inference system, \textit{Knowl.-Based Syst.}, 258 (2022), 110018.

\bibitem{KOKR1} M. Kohler, A. Krzyzak, Over-parametrized deep neural networks do not generalize well, arXiv preprint \url{arXiv:1912.03925}, 2019.

\bibitem{krivoshein2022wavelet} A. Krivoshein, M. Skopina, Wavelet approximation in Orlicz spaces, \textit{J. Math. Anal. Appl.}, 516 (2022), 1, 126473.

\bibitem{LMP24} Y. Li, S. Lu, P. Mathé, S. V. Pereverzev, Two-layer networks with the ReLU k activation function: Barron spaces and derivative approximation, \textit{Numer. Math.}, 156 (2024), 1, 319--344.

\bibitem{LTY19} B. Li, S. Tang, H. Yu, Better approximations of high dimensional smooth functions by deep neural networks with rectified power units, arXiv preprint \url{arXiv:1903.05858}, 2019.

\bibitem{ML2014} W. A. Majewski, L. E. Labuschagne, On applications of Orlicz spaces to statistical physics, \textit{Ann. Henri Poincaré}, 15 (2014), 6, 1197--1221.

\bibitem{MR08} M. Mihăilescu, V. Rădulescu, Neumann problems associated to nonhomogeneous differential operators in Orlicz-Sobolev spaces, \textit{Ann. Inst. Fourier}, 58 (2008), 6, 2087--2111.

\bibitem{Musielak83} J. Musielak, Orlicz spaces and modular spaces, Springer, 2006, vol. 1034.

\bibitem{Musielak87} H. Musielak, Jackson type inequalities for averaged moduli of smoothness, \textit{Comment. Math.}, 27 (1987), 1, Polskie Towarzystwo Matematyczne.

\bibitem{mo59} J. Musielak, W. Orlicz, On modular spaces, \textit{Studia Math.}, 18 (1959), 1, 49--65.

\bibitem{Nikolova2004}
M. Nikolova, A variational approach to remove outliers and impulse noise, {\em J. Math. Imaging Vision}, 20 (2004), 99--120.

\bibitem{OSZ21} J. A. Opschoor, Ch. Schwab, J. Zech, Exponential ReLU DNN expression of holomorphic maps in high dimension, \textit{Constructive Approx.}, 55 (2022), 1, 537--582.

\bibitem{Plesniak1985}
W. Plesniak, Leja's type polynomial condition and polynomial approximation in Orlicz spaces, {\em Ann. Pol. Math.}, 46(1) (1985), 265--275.

\bibitem{plonka2023spline} G. Plonka, Y. Riebe, Y. Kolomoitsev, Spline representation and redundancies of one-dimensional ReLU neural network models, \textit{Anal. Appl.}, 21 (2023), 01, 127--163.

\bibitem{qian2022rates} Y. Qian, D. Yu, Rates of approximation by neural network interpolation operators, \textit{Appl. Math. Comput.}, 418 (2022), 126781.

\bibitem{rao1991theory} M. M. Rao, Z. D. Ren, Theory of Orlicz spaces, Monographs and Textbooks in Pure and Applied Mathematics, 146. Marcel Dekker, New York, 1991.

\bibitem{R2000} M. Ruzicka, Electrorheological fluids: modeling and mathematical theory, Springer, 2007.

\bibitem{SP88} P. Sendov, The Averaged Moduli of Smoothness, Pure and Applied Mathematics, Wiley, 1988.

\bibitem{stein69} E. M. Stein, Note on the class LlogL, \textit{Studia Math.}, 1969.

\bibitem{stein1970singular} E. M. Stein, Singular integrals and differentiability properties of functions, Princeton University Press, 1970.

\bibitem{Y21} D. Yarotsky, Universal approximations of invariant maps by neural networks, \textit{Constructive Approx.}, 55 (2022), 1, 407--474.

\bibitem{ZH1} D.-X. Zhou, Universality of deep convolutional neural networks, \textit{Appl. Comput. Harmonic Anal.}, 48 (2020), 2, 787--794.

\bibitem{Zhikov} V. V. Zhikov, Averaging of functionals of the calculus of variations and elasticity theory, \textit{Math. USSR-Izv.}, 29 (1987), 1, 33.

\bibitem{GNECCObook} R. Zoppoli, M. Sanguineti, G. Gnecco, T. Parisini, Neural approximations for optimal control and decision, Springer, 2020.

\end{thebibliography}
\end{document}